\documentclass[10pt, reqno]{article}

\usepackage[mathscr]{eucal}
\usepackage[all]{xy}
\usepackage{mathrsfs}
\usepackage{xypic}
\usepackage{amsfonts}
\usepackage{amsmath}
\usepackage{amsthm}
\usepackage{amssymb}
\usepackage{latexsym}
\usepackage{eufrak}

\theoremstyle{plain}
\newtheorem{theorem}[equation]{Theorem}
\newtheorem{proposition}[equation]{Proposition}
\newtheorem{lemma}[equation]{Lemma}
\newtheorem{corollary}[equation]{Corollary}

\theoremstyle{definition}
\newtheorem{definition}[equation]{Definition}
\newtheorem{conjecture}[equation]{Conjecture}

\theoremstyle{remark}

\newtheorem{remark}[equation]{Remark}
\newtheorem{claim}[equation]{Claim}

\makeatletter
\renewcommand{\subsection}{\@startsection{subsection}{2}{0pt}{-3ex
plus -1ex minus -0.2ex}{-2mm plus -0pt minus
-2pt}{\normalfont\bfseries}} \makeatother

\numberwithin{equation}{subsection}

\newcommand{\Lmod}[1]{#1\text{-}{\mathsf{mod}}}

\newcommand{\Lmof}[1]{#1\text{-}{\mathsf{mof}}}
\newcommand{\step}[1]{\noindent\vskip 2pt{{\text{\sc{Step }}#1}\en}}

\DeclareMathOperator{\Tor}{\mathrm{Tor}}
\DeclareMathOperator{\Ext}{\mathrm{Ext}}
\DeclareMathOperator{\Hom}{\mathrm{Hom}}
\DeclareMathOperator{\End}{\mathrm{End}}

\DeclareMathOperator{\Ker}{\mathrm{Ker}}

\DeclareMathOperator{\gr}{\mathrm{gr}}

\DeclareMathOperator{\Sym}{\mathrm{Sym}}
\DeclareMathOperator{\supp}{\mathrm{Supp}}
\DeclareMathOperator{\Ad}{\mathrm{Ad}}
\DeclareMathOperator{\Tr}{\mathrm{Tr}}
\DeclareMathOperator{\Lie}{\mathrm{Lie}}
\DeclareMathOperator{\CC}{\mathrm{Ch}}
\DeclareMathOperator{\ad}{\mathrm{ad}}

\newcommand{\dis}{\displaystyle}
\newcommand{\erem}{\hfill$\lozenge$\end{remark}}
\newcommand{\beq}{\begin{equation}\label}
\newcommand{\eeq}{\end{equation}}

\newcommand{\f}[1]{\mathfrak{#1}}
\newcommand{\scr}[1]{\mathscr{#1}}

\DeclareMathOperator{\Spec}{\mathrm{Spec}}
\DeclareMathOperator{\pr}{pr}
\newcommand{\iso}{{\,\stackrel{_\sim}{\longrightarrow}\,}}
\def\eu{{\boldsymbol{1}}}
\def\euu{{\mathsf{eu}}}

\def\ccirc{{{}_{^{\,^\circ}}}}

\newcommand{\GL}{GL}

\newcommand{\sminus}{\smallsetminus}
\renewcommand{\mid}{\enspace\big|\enspace}
\renewcommand{\o}{{\otimes}}
\newcommand{\too}{\,\,\longrightarrow\,\,}
\newcommand{\tooo}{{\;{-\!\!\!-\!\!\!-\!\!\!-\!\!\!\longrightarrow}\;}}

\newcommand{\longeq}{{{=\!\!\!=\!\!\!=\!\!\!=\!\!\!=\!\!\!=\!\!\!=}\,}}

\newcommand{\mto}{\longmapsto}
\newcommand{\map}{\longrightarrow}
\newcommand{\onto}{\,\twoheadrightarrow\,}
\newcommand{\cd}{\!\cdot\!}
\newcommand{\dash}{\mbox{-}}
\newcommand{\pb}{\noindent$\bullet\quad$\parbox[t]{140mm}}

\def\hp{\hphantom{x}}

\renewcommand{\lll}{{\scr L}}
\newcommand{\vi}{$\en{\sf {(i)}}\;$}
\newcommand{\vii}{$\;{\sf {(ii)}}\;$}
\newcommand{\viii}{${\sf {(iii)}}\;$}

\newcommand{\g}{{\mathfrak{g}}}
\newcommand{\fp}{{\mathfrak{p}}}
\newcommand{\bpi}{{\boldsymbol{\pi}}}
\newcommand{\SF}{{\mathsf{F}}}

\newcommand{\triv}{{\mathsf{triv}}}
\newcommand{\A}{{\mathfrak{A}}}

\newcommand{\II}{{\mathbb{I}}}
\newcommand{\I}{{\mathcal{I}}}

\newcommand{\M}{{\mathscr{M}}}
\newcommand{\red}{{\operatorname{red}}}

\newcommand{\inv}{^{-1}}
\newcommand{\oo}{{\mathcal O}}

\renewcommand{\part}{{\f P}}

\newcommand{\slv}{{\mathfrak{s}\mathfrak{l}}(V)}
\newcommand{\sv}{{\mathfrak{s}\mathfrak{l}}}
\newcommand{\mnil}{{\mathscr{M}}_{\mathsf{nil}}}
\newcommand{\mn}{{\mathscr{N}}}
\newcommand{\md}{{\scr C_c}}
\newcommand{\dx}{\D(\X,c)}

\newcommand{\ff}{{\mathbb{F}}}
\newcommand{\BH}{{\mathbb{H}}}
\newcommand{\bv}{{\mathbf{v}}}
\newcommand{\x}{{\mathbf{x}}}
\newcommand{\y}{{\mathbf{y}}}
\renewcommand{\c}{\C[\x,\y]}
\newcommand{\bbf}{{\mathbf{f}}}
\newcommand{\ffs}{{\mathbf{s}}}

\def\C{{\mathbb{C}}}
\def\V{{\f G}}
\def\BN{{\mathsf{N}}}
\def\BT{{\mathsf{T}}}

\def\Q{{\mathsf{Q}}}
\def\pp{{\mathsf{P}}}
\def\zz{\mathcal{Z}}

\def\gln{{\mathfrak{g}\mathfrak{l}}_n(\C)}
\def\gl{{\mathfrak{g}\mathfrak{l}}}

\def\fz{\f z}
\def\fzr{\f z^\circ}
\newcommand{\hh}{{\mathsf{H}}}

\newcommand{\eps}{\epsilon}
\newcommand{\beps}{{\boldsymbol{\varepsilon}}}
\newcommand{\bg}{{\mathfrak{a}}}

\newcommand{\ehe}{{\mathsf{e}\mathsf{H}_c\mathsf{e}}}
\newcommand{\e}{{\mathsf{e}}}
\newcommand{\h}{{\mathfrak{h}}}
\newcommand{\sset}{\subset}
\newcommand{\hreg}{{\mathfrak{h}}^{\operatorname{reg}}}
\newcommand{\bz}{{\mathsf{Z}}}
\newcommand{\opp}{{\operatorname{op}}}
\newcommand{\en}{\enspace}
\newcommand{\D}{{\scr D}}
\newcommand{\bd}{{\mathbf{d}}}
\renewcommand{\P}{{\mathbb{P}}}
\newcommand{\X}{{\f X}}
\newcommand{\xreg}{{\f X}^{\text{reg}}}

\newcommand{\Th}{{\Theta}}
\newcommand{\into}{{\,\hookrightarrow\,}}

\newcommand{\G}{\Gamma}
\newcommand{\Om}{\Omega}
\newcommand{\OO}{{\mathsf{O}}}
\newcommand{\OOO}{{\overline{\mathsf{O}}_1}}
\newcommand{\ind}{{\f I}}
\newcommand{\Ug}{{\mathcal{U}\f g}}
\newcommand{\U}{{\mathcal{U}}}

\newcommand{\Z}{{\mathbb{Z}}}
\newcommand{\La}{\Lambda}
\newcommand{\gc}{{\g_c}}
\newcommand{\ham}{{\mathbb{H}}}
\newcommand{\Id}{\operatorname{Id}}
\newcommand{\Hilb}{{\operatorname{Hilb}^n{\mathbb{\C}}^2}}
\newcommand{\fl}{{\f l}}
\newcommand{\bu}{{\boldsymbol{\mathcal{U}}}}
\newcommand{\rel}{{\X^{\operatorname{relevant}}}}
\newcommand{\scal}{{\scr S}}
\newcommand{\eer}{{{\mathfrak{G}}^{\text{reg}}}}

\newcommand{\N}{{\mathbb{N}}}
\newcommand{\al}{{\alpha}}
\newcommand{\be}{{\beta}}
\newcommand{\la}{{\lambda}}
\newcommand{\Mat}{{\mathrm{Mat}}}
\newcommand{\QQ}{{\overline{Q}}}
\newcommand{\Rep}{{\mathrm{Rep}}}
\renewcommand{\SS}{{\overline{S}}}
\newcommand{\aug}{\langle I \rangle}
\begin{document}
\centerline{\Large{\textbf{\Large{Almost-commuting variety, $\D$-modules,
and Cherednik Algebras}}}}

\bigskip

\centerline{\sc Wee Liang Gan and
Victor Ginzburg}
\medskip

\begin{abstract} We study a scheme $\M$ closely related to the
set of  pairs
of $n\times n$-matrices with rank 1 commutator. We show that 
 $\M$  is a reduced complete intersection with
$n+1$ irreducible components, which we
describe.

There is a distinguished  Lagrangian subvariety
 $\mnil\sset \M$. We introduce a category, $\scr C$,
of $\D$-modules
whose {\em characteristic variety} is contained in 
$\mnil$. Simple objects of that category are
analogous to Lusztig's {\em character sheaves}.
We construct an exact functor of {\em Quantum Hamiltonian
reduction} from our category $\scr C$ to the
category $\mathcal O$ for type $\mathbf{A}$ rational
Cherednik algebra. 
\end{abstract}

\centerline{\sf Table of Contents}
\vskip -1mm

$\hspace{30mm}$ {\footnotesize \parbox[t]{115mm}{
\hp${}_{}$\!
\hp\!1.{ $\;\,$} {\tt Introduction.} \newline
\hp2.{ $\;\,$} {\tt The geometry of $\M$.} \newline
\hp3.{ $\;\,$} {\tt Generalization to quiver moment  maps.} \newline
\hp4.{ $\;\,$} {\tt A Lagrangian variety.} \newline
\hp5.{ $\;\,$} {\tt A category of holonomic $\D$-modules.}\newline
\hp6.{ $\;\,$} {\tt Cherednik algebra and Hamiltonian reduction.}\newline
\hp7.{ $\;\,$} {\tt The functor of Hamiltonian reduction.}\newline
\hp8.{ $\;\,$} {\tt Appendix by {V. Ginzburg:}\en{A remark on a theorem
of {M. Haiman.}}}
}
}

\section{Introduction}
\subsection{}\label{int1}
Let $V:=\C^n$ and let  $\g := \End(V)=\gln$ 
be the Lie algebra of $n\times n$-matrices.
We will write 
elements of $V$ as column vectors, and elements of
$V^*$ as row vectors.
We consider the following
affine closed subscheme in the vector space
$\g\times\g\times V\times V^*$:
\begin{equation}\label{M}
 \M := \{ (X,Y,i,j)\in \g\times\g\times
V\times V^* \mid [X,Y]+ij=0\}. 
\end{equation}
More precisely, let 
$\C[\g\times\g\times V\times V^*]=\C[X,Y,i,j]$ denote the polynomial
algebra, and let $J\subset \C[X,Y,i,j]$ be the ideal
generated by the $n^2$ entries of the matrix $[X,Y]+ij$.
Then, by definition, we have $\M=\Spec\C[X,Y,i,j]/J,$
a not necessarily reduced affine scheme.

All the results of this paper are based on 
Theorem \ref{t1} below that describes the structure of the
scheme $\M$.
To formulate the Theorem,
for each integer $k\in\{0,1,\ldots,n\}$, let
$$ \M'_k := \bigg\{ (X,Y,i,j)\in\M \left |
\begin{array}{c}
\textrm{$Y$ has pairwise distinct eigenvalues},\\
\dim(\C[X,Y]i)=n-k,\quad \dim(j\C[X,Y])=k \end{array}\right\} $$
and let $\M_k$ be the closure of $\M'_k$ in $\M$.


\begin{theorem} \label{t1}
\vi $\M$ is a complete intersection in $\g\times\g\times
V\times V^*$.

\vii The irreducible components of $\M$ are $\M_0, \ldots, \M_n$.

\viii $\M$ is reduced and equidimensional; we have $\dim\M=n^2+2n$.
\end{theorem}
As we have learned after completion of the paper,
a description of irreducible components of a scheme closely
related to $\M$ was found earlier by   Neubauer \cite{Ne}.

We may identify $\g$ with its dual, $\g^*,$
and write $\g\times\g\times V\times V^*\cong T^*(\g\times V).$
The cotangent bundle $T^*(\g\times V)$ comes equipped
with the standard symplectic structure.
A (possibly singular or non-reduced) subscheme $Z\sset T^*(\g\times V)$
is said to be {\em Lagrangian} if the generic locus
of each irreducible component of $Z_{\operatorname{red}}$,
the scheme $Z$ taken with reduced structure,  is a Lagrangian
subvariety in $T^*(\g\times V)$.

We 
introduce the following closed subset of $\M$, to be given  a scheme
 structure later, in \S\ref{mnil_sec},
\begin{equation}\label{nil}
\mnil:=\{(X,Y,i,j)\in \g\times\g\times V\times V^*
\mid [X,Y]+ij=0\enspace\&\enspace Y\enspace \text{is nilpotent}\}.
\end{equation}

\begin{theorem} \label{nilt}
The scheme $\mnil$ is a
{\sf{Lagrangian}}  (not necessarily reduced)
complete intersection in $T^*(\g\times V).$
\end{theorem}

\subsection{Applications to the commuting variety.}\label{int2}
Write $\zz=\{(X,Y)\in\g\times\g\mid [X,Y]=0\}$
for the {\em commuting variety}.  Again,
we regard $\zz$ as a (not necessarily reduced)
affine subscheme
$\zz:=\Spec\C[X,Y]/I,$ where
$\C[X,Y]=\C[\g\times\g]$ stands for the
polynomial algebra, and $I\subset \C[X,Y]$ stands for the
ideal
generated by the $n^2$ entries of the matrix $[X,Y]$.
The scheme $\zz$ is known to be irreducible (cf. \cite{Ri}).

Further, write $G:=\GL(V)$. The group
$G$ acts diagonally on $\g\times\g$ via
$ g\cdot(X,Y) = (gXg^{-1}, gYg^{-1})$.
The induced $G$-action on $\C[X,Y]$ by
algebra automorphisms clearly preserves the
ideal $I$.

It is a well known and
 long standing open question whether or not the scheme
$\zz$ is reduced,
i.e.,  whether or not $\sqrt{I}=I$.
We cannot resolve this question. However, 
during discussions with Pavel Etingof we have realized that
Theorem \ref{t1} combined with some elementary results from
\cite{EG} implies the following

\begin{theorem} \label{t2}
One has: $I^G = (\sqrt{I})^G$.
\end{theorem}

\begin{remark} Write $\zz_\red$ for the scheme $\zz$ taken with 
reduced scheme structure. It follows from Theorem \ref{t1}(ii)-(iii) that
the map $(X,Y,i,j)\mapsto (X,Y,i)$ gives
 an isomorphism of algebraic varieties $\M_0\simeq \zz_\red\times
V.$
\end{remark}

\subsection{Cherdnik algebras and quantum Hamiltonian reduction.}
\label{hhham}
Our interest in the geometry discussed in \S\ref{int1} comes
from the theory of {\em rational Cherednik algebras},
 an important class of associative
algebras introduced 
in \cite{EG}. Below, we will
only consider  rational Cherednik algebras of type $\mathbf{A}$.
Specifically,  let $\h=\C^n$ be  the tutological 
$n$-dimensional permutation representation of the Symmetric group
$S_n$. The  corresponding Cherednik algebra
is generated by a copy of the vector space $\h$,
a copy of the dual space $\h^*,$ and also by the
elements $w\in S_n$. These
generators are subject to the defining relations
\eqref{defrel},
in \S\ref{reminder}. The resulting algebra $\hh_c$
is referred to as the  rational Cherednik algebra of $\gl_n$-type
with parameter $c\in\C$.

Write $\e:=\frac{1}{n!}\sum_{\sigma\in 
S_n}\,\sigma\in \C[S_n]\subset\hh_c\,$ 
for the symmetrizer idempotent. 
The subalgebra $\ehe\sset \hh_c$ is called 
the {\em spherical subalgebra} of $\hh_c$.
It has been argued  in \cite{EG} that
the spherical subalgebra may be viewed as a quantization
of  so-called {\em Calogero-Moser space}. Specifically, given  $c\in\C^*$,
let  $\OO_c:=\Ad G(\chi_c)\sset\g$, be the semisimple conjugacy class of the matrix
$\chi_c:=c(\Id -n\cdot p),$ where $p=\operatorname{diag}(0,0,\ldots,0, 1)$ is the projector
on the line $\ell$ spanned by the last coordinate vector.
Then, by an old result due to Kazhdan-Kostant-Sternberg
\cite{KKS}, the
 Calogero-Moser space (with parameter $c$) may be interpreted as a
{\em classical} Hamiltonian reduction of
the symplectic vector space $\g\times\g^*$ over  $\OO_c$, viewed
as a coadjoint orbit in $\g^*\cong\g$.

Now, the Poisson algebra $\C[\g\times\g^*]$ has a natural
quantization $\D(\g)$, the  algebra of polynomial differential operators 
on $\g$.
Accordingly, one of the main results of \cite{EG} says that, for generic values
of $c$, the algebra  $\ehe$ may be constructed as a {\em quantum}
 Hamiltonian reduction of $\D(\g)$ with respect to
 the primitive ideal in the enveloping algebra of $\g$ that
`quantizes' the coadjoint orbit $\OO_c$.
It has been conjectured that a similar result should hold, in effect,
for all values of $c$.
This conjecture is proved in the present paper in full generality,
 see Theorem \ref{ker}.

The main novelty of our
present approach, as compared to that of \cite{EG},
 is in
replacing the  coadjoint orbit $\OO_c=\Ad G(\chi_c)$ by its standard
{\em polarization}, an affine Lagrangian subspace
of the form $\chi_c+\fp^\perp,$ where 
$\fp\sset \g$ is the maximal parabolic formed by the
matrices which preserve the line $\ell\sset  V$.
It is a well known heuristic
general principle of representation theory
that, usually,
performing Hamiltonian reduction (either
classical or quantum) with respect to
the group $G$ and its coadjoint orbit $\OO_c$ should be equivalent
to performing Hamiltonian reduction with respect to
the subgroup $P\sset G$, corresponding to the
Lie algebra $\fp$, and the one-point 
coadjoint $P$-orbit $\{\chi_c|_{_\fp}\}\sset\fp^*$.
Applying this  heuristic principle in our situation,
 we are thus
led to consider  quantum Hamiltonian reduction
of  the algebra $\D(\g)$ with respect to
the group $P$ and its character $\chi_c$. Observe that the group
$P$ was defined as the isotropy group of $\ell$, a point in the
projective space $\P=\P(V)$.
Therefore, a routine argument
shows that performing Hamiltonian reduction 
of 
the algebra  $\D(\g)$  with respect to $P$ is the same thing
as  performing Hamiltonian reduction of 
 $\D(\g\times\P)$, a  larger algebra,
with respect to the group $G$ acting diagonally on 
$\g\times\P$.
It is this latter reduction that we are doing in the present paper. 

Let  $\D(\g\times\P,c)$ be the algebra of $c$-{\em twisted}
 differential operators
on $\g\times\P$, and let $\gc$ denote the image of the 
Lie algebra $\slv\sset\g=\gl(V)$ in  $\D(\g\times\P,c)$,
 see \S\S\ref{tdo},\ref{pf} for more details and unexplained notation.
Our main result stated below  gives a construction of the
spherical subalgebra $\ehe$ in terms of 
quantum Hamiltonian reduction
of the algebra $\D(\g\times\P,c)$.

\begin{theorem}\label{iso} For any $c\in\C$, let
$\hh_c$ be the rational Cherednik algebra of $\gl_n$-type with parameter
$c$. Then, there is 
a filtered algebra isomorphism
\[\Phi_c: \
\bigl(\D(\g\times\P,c)/\D(\g\times\P,c)\cdot{\gc}\bigr)^{\ad\gc}\iso\ehe
\]
such that 
\begin{equation}\label{psi}
\Phi_c(\C[\g]^{\Ad G})=\C[\h]^{S_n}\sset\ehe,\quad\text{and}\quad
\Phi_c(\bz)=(\Sym\h)^{S_n}\sset\ehe.
\end{equation}

Moreover, the associated graded map gives a graded algebra
isomorphism 
\[\gr\Phi_c:\
\gr\bigl(\D(\g\times\P,c)/\D(\g\times\P,c)\cdot{\gc}\bigr)^{\ad\gc}
\iso  \gr(\ehe).\]
\end{theorem}

Theorem  \ref{iso} is a strengthening of
\cite[Section~7]{EG}, esp. Corollary 7.4, (cf. also \cite{BFG}, Theorem
7.2.4(i)).
In \cite{EG},
the homomorphism $\Phi_c$ has been called {\em Harish-Chandra isomorphism}
for Cherednik algebras.

The  Cherednik algebra of $\gl_n$-type  contains,
 as a subalgebra, the  Cherednik algebra of
${\mathfrak{s}\mathfrak{l}}_n$-type.
The latter is obtained by replacing,
in the definition of $\hh_c$,
the $n$-dimensional $S_n$-representation
$\h=\C^n$ by its $(n-1)$-dimensional irreducible
subrepresentation.
Our proof of Theorem \ref{iso} applies also to
the  rational Cherednik algebra
of ${\mathfrak{s}\mathfrak{l}}_n$-type. In that case, the sourse of
the homomorphism $\Phi_c$ should be 
the algebra
$\bigl(\D({\mathfrak{s}\mathfrak{l}}_n\times \P,c)/\D({\mathfrak{s}\mathfrak{l}}_n\times
\P,c)\cdot{\gc}\bigr)^{\ad\gc},$
a Hamiltonian reduction of the algebra
$\D({\mathfrak{s}\mathfrak{l}}_n\times \P,c)$.

\subsection{$\G$-analogues.} 
Theorem \ref{t1} has a 
$\G$-equivariant generalization,
where
 $\G\subset SL_2(\C)$ is
 a finite subgroup.

Fix an integer $n\geq 1$ and
 let ${\mathbf{R}}=R^{\oplus n}$
be the direct sum of $n$ copies of the left regular  $\G$-representation.
Further, let $x,y$ be the standard basis in $\C^2$, the tautological
2-dimensional $\G$-representation.
Using this basis, any linear map
$F:{\mathbf{R}}\otimes_{_\C}\C^2\to {\mathbf{R}}$
may be identified with a pair of linear maps
$X:=F(-\otimes x),\,Y:=F(-\otimes y):\, {\mathbf{R}}\to {\mathbf{R}}.$
Also, given a vector $i\in {\mathbf{R}}$ and
 a covector  $j\in {\mathbf{R}}^*:=\Hom_\C({\mathbf{R}},\C),$
we have $i\otimes j\in {\mathbf{R}}\otimes{\mathbf{R}}^*=\End_{_\C}{\mathbf{R}}.$

As a generalization of \eqref{M},
one introduces the 
following affine scheme, see \cite{Na1} and also \cite[formula (1.11)]{EG}:
$$
\M(\G,n) :=
\big\{(X,Y,i,j)\in \bigl(\Hom_{_\C}({\mathbf{R}}\otimes_{_\C}L,{\mathbf{R}})
\times {\mathbf{R}}\times
{\mathbf{R}}^*\bigr)^\G\enspace
\Big|\enspace
[X,Y] + i\otimes j = 0\big\}.
$$

The following theorem reduces, in the special
case $\G=\{1\}$, to parts (i) and (iii) of Theorem~\ref{t1}.

\begin{theorem}\label{tquiver} 
The scheme $\M(\G,n)$ is reduced and has $n+1$ irreducible components.
Moreover, it
 is a complete intersection in 
$\bigl(\Hom_{_\C}({\mathbf{R}}\otimes_{_\C}L,{\mathbf{R}})
\oplus {\mathbf{R}}\oplus
{\mathbf{R}}^*\bigr)^\G$.
\end{theorem}

We shall  describe the irreducible components of
$\M(\G,n)$ in Section 3.
The proof of Theorem \ref{tquiver} 
is similar to the proof of Theorem \ref{t1};
the difference is that it makes use of
results by Crawley-Boevey \cite{CB} which we were able to
avoid in the special case of $\Gamma=\{1\}$.

Theorem \ref{tquiver} plays an important role in the
construction of Harish-Chandra homomorphism for
wreath-products, see \cite{EGGO}, that reduces
in the special case of   $\Gamma=\{1\}$
to  Theorem \ref{iso} above.

For any {\em nontrivial}
finite subgroup  $\G\subset SL_2(\C)$,
G. Lusztig has  considered in  \cite[\S12]{Lu3}
a certain closed 
subset $\mnil^{\text{Lus}}(\G,n)\sset\M(\G,n)$ 
and proved that this subset is  a Lagrangian subscheme.
In the special  case of the trivial group
$\G=\{1\}$, the set  $\mnil^{\text{Lus}}(\G,n)$ is still well defined.
It turns out that this set is properly contained in, but is
{\em not} equal to, our set $\mnil$, see \eqref{nil}.
Roughly speaking, the difference between the two sets is  that,
 in our definition,
only the operator $Y$ in the quadruple $(X,Y,i,j)$
is required to be nilpotent, 
while in Lusztig's  definition both $X$ and $Y$ are
 required to be nilpotent. 

More generally, a Lagrangian scheme $\mnil(\G,n)\sset\M(\G,n),$ that properly
contains   Lusztig's Lagrangian
scheme  $\mnil^{\text{Lus}}(\G,n)$ and which is analogous to our scheme    \eqref{nil},
may be defined in the case of  any
 {\em cyclic} (in particular, trivial)
group $\G$, that is, in the case of quivers of type $\widetilde{\mathbf{A}}$.
On the other hand, there seems to be no
analogue of such a scheme in the
case of  finite subgroups $\G\sset SL_2(\C)$
of non-cyclic type.

\begin{remark} In the quiver language,
 the  case of the trivial group $\G=\{1\}$   corresponds to the quiver
with a single vertex and a single edge-loop at that vertex.
This case 
 does  {\em not} fall in the setting of  \cite{Lu3}, since
Lusztig  excludes the case of  quiver with edge-loops.
\hfill$\lozenge$\end{remark}

\subsection{} Here are  more details about the contents
of the paper.

In \S2, we prove Theorems \ref{t1}, \ref{nilt},
and Theorem \ref{t2}. In \S3, we prove Theorem \ref{tquiver}.
In \S4 we study a certain
Lagrangian variety, $\La$, closely related to the variety $\mnil$.

In $\S5$, for each $c\in \C$,
we introduce a category  $\scr C_c$ of 
holonomic  $\D(\g\times\P,c)$
whose characteristic variety is contained in $\La$.
 Associated with $\La$, there is a natural stratification of
the space $\g\times \P(V)$, and $\D$-modules from
category  $\scr C_c$ are smooth along the strata
of that  stratification.
Simple objects of  category  $\scr C_c$ are
analogous to Lusztig's {\em character sheaves}, see \cite{Lu2}.
More results about category  $\scr C_c$ will be given in \cite{Gi3}.

In Section 6, we remind the definition of
 rational
Cherednik algebra of type $\mathbf{A}_{n-1}$,
 and  prove Theorem \ref{iso}.
In \S7, we introduce and study a natural
 {\em Hamiltonian
reduction functor}
 $\ham: \scr C_c\map \mathcal O(\ehe)$, 
an exact functor from category $\scr C_c$ to the
category $\mathcal O$ for the spherical subalgebra $\ehe$, as defined
e.g., in \cite{BEG}. In the Appendix (\S8), the geometry of the variety
$\M$ is applied to deduce a result of M. Haiman on
powers of the space of alternating (with respect to the
 $S_n$-diagonal action) polynomials
in $2n$ variables.

\subsection{Acknowledgements.}{\small{ The  second author would like
to thank   Alexander Kuznetsov from whom he learned  some years ago
that the scheme $\M$ should have
$n+1$ irreducible components. 
Both authors are  grateful to Pavel Etingof
for very useful discussions which have led us to
Theorem \ref{t2}.
Thanks are also due to Rupert Yu for pointing out a mistake in an
earlier version of the paper,  to Michael Finkelberg,
and Iain Gordon for interesting
correspondence and  comments, 
and to Robert Guralnick for bringing references
\cite{Ne} and \cite{Gu} to our attention.
 The work of W.L.G. was partially supported by the
NSF grant  DMS-0401509. The work of V.G. was partially supported by the NSF grant
DMS-0303465 and CRDF grant RM1-2545-MO-03. 
}}

\section{The geometry of $\M$}\label{geometry}

This section is devoted to the proofs of  Theorem \ref{t1} and Theorem \ref{t2}.

\subsection{Linear algebra.}\label{proofs1}
Given  $X\in\g$, let $G^X:=\{g\in  \GL(V)\mid gX=Xg\}$ be the
centralizer of $X$ in $G=\GL(V)$.

\begin{lemma} \label{ee} For any $X\in\g$, the group
$G^X$ acts on $V$ with finitely many orbits.
\end{lemma}

\proof We may assume that $X$ is in Jordan normal form.
Suppose first that there is only one Jordan block,
and the corresponding eigenvalue of $X$ is $\lambda$.
In this case, the non-zero orbits of $G^X$ are the sets
$$\{ v\in V \mid (X-\lambda)^k v \neq 0\enspace\&
\enspace
 (X-\lambda)^{k+1} v =0 \}
\quad \mbox{ where } k=0,1,\ldots,n-1.$$

Let now $X$ have  several Jordan blocks. We write a direct sum decomposition
$V=\oplus_s V_s$ according to the block decomposition of $X$,
and let $\overline{G}^X:=G^X\,\bigcap\, \bigl(\prod_s \GL(V_s)\bigr)$
be the part of the group $G^X$ that respects the
 direct sum decomposition. We deduce from the above
that the group $\overline{G}^X$, hence the larger group $G^X$,
has finitely many orbits in $V=\bigoplus_s V_s$.
\endproof

Let  $\mn\sset\g$ be the nil-cone formed by nilpotent matrices,
and  let $G$ act diagonally on  $\mn\times V$.

\begin{corollary}\label{mnfinite} The set  $\mn\times V$ is a finite 
union of $G$-diagonal orbits.
\end{corollary}
\proof 
The nil-cone $\mn$ is a finite union of $G$-orbits,
hence we have a finite partition
$\mn\times V=\bigsqcup\,\OO\times V$,
where $\OO$ runs over the  $G$-orbits in $\mn$.
Now, each set $\OO\times V$ is clearly stable
under the $G$-diagonal action.
Furthermore, $G$-diagonal orbits in  $\OO\times V$
are in one-to-one correspondence with
$G^X$-orbits in $V$, where $X$ is  some fixed
element of $\OO$. Thus, we are done
by  Lemma~\ref{ee}.
\endproof

Fix a quadruple $(X,Y,i,j)\in \g\times \g\times V\times V^*.$
Denote 
by $\C[X,Y] i$ the subspace of $\C^n$ consisting of vectors
of the form $Ai$, where $A$ is any matrix which can be written
as a noncommutative polynomial in $X$ and $Y$. 
The following lemma was due to Nakajima \cite[Lemma 2.9]{Na}; 
we give an alternative shorter proof.

\begin{lemma} \label{nak}
If $(X,Y,i,j)\in\M$, then $j$ vanishes on $\C[X,Y] i$.
\end{lemma}
\proof
Since the rank of the matrix $[X,Y]$ is at most $1$, we 
can simultaneously conjugate $X$, $Y$ into upper triangular
matrices, cf. \cite{Gu} and also \cite[Lemma 12.7]{EG}. Hence, we assume
without loss of generality that $X$, $Y$ are upper triangular
matrices. In this case, for any $A\in\C[X,Y]$, we have
\beq{ji}
jAi = \Tr(Aij) = -\Tr(A[X,Y]) = 0
\eeq
since $A$ is upper triangular and $[X,Y]$ is strictly
upper triangular.
\endproof

\subsection{The moment map.}
We will  identify $\g^*$ with $\g$ via the pairing
$\g \otimes \g \too \C,\, X\otimes Y \mapsto \Tr(XY)$. 
Let $\V := \g\times V$, and view it as a $G$-variety with respect to
the $G$-diagonal action.

The  induced $G$-action on
$T^*\V \simeq \g\times\g\times V\times V^*$ is given by the formula
$$ g\cdot(X,Y,i,j) = (gXg^{-1}, gYg^{-1}, gi, jg^{-1})
\quad \mbox{where $g\in G$}.$$
This $G$-action on $T^*\V$ is hamiltonian and the
corresponding  moment map is 
\begin{equation}\label{moment}
\mu: T^*\V=\g\times\g\times V\times V^* \too \g^*\simeq\g,\quad
\quad (X,Y,i,j)\longmapsto [X,Y]+ij.
\end{equation}

We see that the scheme $\M$, see \eqref{M},
may (and will) be identified with $\mu^{-1}(0)\sset T^*\V$,
 the scheme-theoretic zero fiber of
the moment map \eqref{moment}.

 Given a conjugacy class $\OO$ in $\g$,
we put $\M(\OO):=\{(X,Y,i,j)\in\M\,|\, Y\in\OO\},$
viewed as a (not necessarily closed)
reduced scheme.

\begin{proposition}\label{lagrange} $\M(\OO)$
is a Lagrangian subscheme in~$T^*\V$, for any conjugacy class
$\OO\sset\g$.
\end{proposition}

\proof In general, assume first that
 $\V$ is an arbitrary smooth $G$-variety.
Write $\scal$ for the set (possibly infinite) of all
$G$-orbits in $\V$ and, given a $G$-orbit
$S\in\scal$, let $T_S^*\V\sset T^*\V$ denote the conormal
bundle to $S$. Then, the natural $G$-action
on $T^*\V$ is Hamiltonian with moment map
$\mu:T^*\V\to\g^*$, and it is well known that
\begin{equation}\label{lagr}\mu^{-1}(0)=\bigcup_{S\in\scal} T_S^*\V.
\end{equation}

Now, return to our case $\V=\g\times V$ and
write $\pr: \g\times V^*\to\g$ for the first
projection. We see from \eqref{lagr} that proving
the Proposition amounts to showing
that, for any
 conjugacy class $\OO$ in $\g$,
the set $\pr^{-1}(\OO)$ is a finite union of $G$-orbits.
It is clear that this last statement is equivalent to 
Lemma \ref{ee}.
\endproof

\subsection{A flat morphism.}
Write $\C^{(n)}$ for the set of unordered $n$-tuples of complex
numbers. Let 
\begin{equation}\label{spec}
\bpi:\ T^*\V=\g\times\g\times V\times V^*\too \C^{(n)},\quad (X,Y,i,j)\mto\Spec Y
\end{equation}
 be the map that sends
$(X,Y,i,j)$ to the unordered $n$-tuple $\Spec Y$ of eigenvalues of $Y$,
counted with multiplicities. 

\begin{proposition}\label{fib_dim} 
The morphism  $\mu\times\bpi: T^*\V\to \g\times \C^{(n)}$
is flat. In particular,  all nonempty  (scheme-theoretic)
fibers of  this morphism are
equidimensional,
of dimension $n^2+n.$
\end{proposition}

\begin{proof}
For any unordered $n$-tuple
$y=(y_1,\ldots,y_n)\in\C^{(n)},$ the set of all the matrices
$Y\in\g$ such that $\Spec Y=y$
is clearly a finite union of  conjugacy classes.
Therefore, the zero fiber of the map $\mu\times\bpi$ is
equal, as a set, to
a finite union of Lagrangian subschemes of the form
$\M(\OO)$, as in Proposition \ref{lagrange}.
In particular, the dimension of  the zero fiber fiber
is $\leq \frac{1}{2}\dim T^*\V=n^2+n=
\dim T^*\V-\dim(\g\times\C^{(n)}). $

Next, we define a $\C^\times$-action on each of the varieties
$T^*\V, \, \g,$ and $\C^{(n)}$, as follows.
We let $z\in\C^\times$ act on
$\g\times\g\times V\times V^*,$ resp.,
on $\C^n,$ by scalar multiplication by $z$.
This gives a  $\C^\times$-action on $T^*\V,$ resp., on
 $\C^{(n)}=\C^n/S_n$.
Further, we let $z\in\C^\times$ act on
$\g$  by scalar multiplication by $z^2$.
This gives a  $\C^\times$-action on $\g$ such
that the map $\mu\times\bpi$ becomes a $\C^\times$-equivariant
morphism. Thus, the standard argument based on
the asymptotic cone construction, cf. eg. \cite[\S2.3.9]{CG} or
\cite[ch.I, \S6]{Kr}, shows that the
dimension of any fiber of  the map $\mu\times\bpi$
cannot be greater than the dimension of the
zero fiber. 

Hence, since $\g\times \C^{(n)}$ is smooth,
we conclude
 that the morphism  $\mu\times\pi$
is flat,  cf. \cite[Proposition 6.1.5]{Gr}.
\end{proof}

Composing the flat morphism
$\mu\times\bpi$ with the first projection
$\g\times \C^{(n)}\onto \g$, we deduce
the following special case of  \cite[Theorem 4.4]{CB}.

\begin{corollary} The morphism $\mu$ is flat.\qed
\end{corollary}

Let $\pi=\bpi|_\M$ be the restriction of $\bpi$ to the
closed subscheme $\M\sset T^*\V$.

\begin{corollary} \label{flat} The scheme $\M$ is a
 complete
intersection in $T^*\V$ and $\dim\M=n^2+2n$. 
Furthermore,
$\pi:\M\to \C^{(n)},\, (X,Y,i,j)\mapsto\Spec Y$ is a flat  morphism.
All fibers of this morphism are $(n^2+n)$-dimensional,
 Lagrangian subschemes in $T^*\V$.
\end{corollary}
\begin{proof}
First of all, any quadruple of the form $(0,Y,0,0)$
belongs to $\mu^{-1}(0).$ It follows that all fibers
of the restriction of the map $\bpi$ to $\M$ are nonempty.

Now, flat base change with respect to the
imbedding $\{0\}\times \C^{(n)}\into\g\times \C^{(n)}$
implies that
the sheme $\M=\mu\inv(0)=(\mu\times\bpi)\inv(\{0\}\times \C^{(n)})$
is a  complete
intersection in $T^*\V$ and, moreover,
that the morphism $\pi:\M\to \C^{(n)}$ is flat.
In particular, the dimension of any irreducible component of any fiber of
this morphism equals
$\dim T^*\V-\dim(\g\times \C^{(n)})= n^2+n$.

Further, it is clear that  each fiber of the map
$\pi$ is equal, as a set, to a finite union of  
 Lagrangian subschemes of the form
$\M(\OO)$, see Proposition \ref{lagrange}.
Furthermore, we have proved that 
each  irreducible component of the corresponding scheme-theoretic
fiber has the same dimension as the dimension of $\M(\OO)$.
Thus, any such irreducible component must be the closure of an
irreducible
component of the set of the
form $\M(\OO),$ hence, it is a Lagrangian subscheme.
\end{proof}

\subsection{The scheme $\mnil$.}\label{mnil_sec} First of all, 
recall that  the nil-cone  $\mn\sset \g$ is equal, as a set,
to the zero fiber of
the map $\g\to\C^{(n)},\,Y\mapsto\Spec Y,$ cf. \eqref{spec}.
Thus, we make $\mn$ a scheme by giving it the
scheme structure of the scheme-theoretic zero fiber of
the map $Y\mapsto\Spec Y.$
It is known that this scheme is an irreducible  reduced
scheme of dimension $n^2-n$;
moreover, it is a complete intersection in  $\g$.

We consider the projection $\pr_Y:\M\to\g, (X,Y,i,j)\mapsto Y$.
It is clear that set-theoretically we have $\mnil=(\pr_Y)\inv(\mn),$
 see \eqref{nil}.
We define a scheme structure on $\mnil$ to be the natural one
on  the {\em scheme-theoretic} inverse image of the scheme $\mn$ under
the
morphism $\pr_Y$. 
We do not know whether or not the scheme  $\mnil$ is reduced.

The above discussion shows that $\mnil$ may be identified, as a scheme,
with the scheme-theoretic zero fiber of the morphism
$\pi: \M \to\C^{(n)}.$ Thus, by Corollary \ref{flat},
$\mnil$ is a Lagrangian complete intersection.
This proves  Theorem \ref{nilt}.\qed

\subsection{Generic locus of $\M$.}
Our proof of Theorem \ref{t1}
follows the strategy of \cite[\S1]{Wi}, in which 
Wilson considered the equation $[X,Y]+ \mathtt{Id}=ij$ instead 
of our equation $[X,Y]+ij=0$. Wilson's situation was somewhat simpler
since in his case the corresponding variety $\M$ was
 {\em smooth} and {\em irreducible}.

\begin{lemma}\label{cc}
Let $(X,Y,i,j)\in \M$. Suppose that
the eigenvalues of $Y$ are pairwise distinct,
$\dim \C[X,Y]i \leq n-k$, and $\dim j\C[X,Y]\leq k$. 
Then the $G$-orbit of $(X,Y,i,j)$ contains a representative
such that:
\begin{itemize}
\item[\vi] $Y$ is diagonal, 
say $Y=\mathrm{diag}(y_1,\ldots,y_n)$;
\item[\vii] $i=(0,\ldots,0,\underbrace{1,\ldots, 1}_{n-k''})$ and
$j=(\underbrace{1,\ldots, 1}_{k'},0,\ldots,0)$
for some $k'\leq k\leq k''$;
\item[\viii] $X=(X_{rs})_{1\leq r,s\leq n}$ has the entries 
$$ X_{rs} = \left\{ \begin{array}{ll}
x_r & \textrm{if $r=s$,}\\
\frac{1}{y_r-y_s} & \textrm{if $r> k''$ and $s\leq k'$,}\\
0 & \textrm{else,}
\end{array} \right. $$
for some $x_1, \ldots, x_n$.
\end{itemize}

\noindent
Conversely, the data $(X,Y,i,j)$ defined by \vi,\vii and \viii 
for any choices of $x_1, \ldots, x_n$, $y_1, \ldots, y_n$, $k'$, $k''$ 
(with $y_1, \ldots, y_n$ pairwise distinct, $k'\leq k\leq k''$)
 belongs to $\M$, moreover, we have 
$\dim(\C[X,Y]i)=n-k''$, $\dim(j\C[X,Y])=k'$.
\end{lemma}
\proof (cf. \cite[Proposition 1.10]{Wi}.)
We may assume that $Y=\mathrm{diag}(y_1,\ldots,y_n)$. 
We may assume furthermore that
$j=(\underbrace{1,\ldots, 1}_{k'},0,\ldots,0)$ for some $k'$.
By the Vandermonde determinant, we see that $j\C[Y]$ is the
$k'$ dimensional subspace of $\C^n$ whose last $n-k'$ coordinates
are $0$. Hence, $k'\leq k$, and by Lemma \ref{nak} we may assume 
that $i=(0,\ldots,0,\underbrace{1,\ldots, 1}_{n-k''})$ for
some $k''\geq k'$. By the Vandermonde determinant again,
we see that $\C[Y]i$ is the $n-k''$ dimensional subspace of 
$\C^n$ whose first $k''$ coordinates are $0$, hence $k''\geq k$.
Now, solving the equation $[X,Y]+ij=0$ for $X$ gives (iii).
Note that since $X$ is lower-triangular, we have
$j\C[X,Y]=j\C[Y]$ is $k'$ dimensional, 
and $\C[X,Y]i=\C[Y]i$ is $n-k''$ dimensional.
Thus, we have proved that the
$G$-orbit has a representative satisfying (i), (ii) and (iii),
and we have also proved the last statement of the lemma.
\endproof

\subsection{Irreducible components of $\M$.}
 For each $k=0,\ldots,n,$ we introduce 
the following subset of $\M$:
$$ \M''_k := \bigg\{ (X,Y,i,j)\in\M \left |
\begin{array}{c}
\textrm{$Y$ has pairwise distinct eigenvalues},\\
\dim(\C[X,Y]i)\leq n-k,\quad \dim(j\C[X,Y])\leq k,\\
\mbox{ and }\dim(\C[X,Y]i)+\dim(j\C[X,Y])<n \end{array}\right\} .$$

\begin{lemma} \label{dd}
\vi The dimension of $\M'_k$ is equal to 
$n^2+2n$, and $\M'_k$ is connected. Moreover, both the actions 
of $G$ and $\g$ on $\M'_k$ are free.

\vii The dimension of $\M''_k$ is strictly less than $n^2+2n$.
\end{lemma}
\proof
Consider any $(X,Y,i,j)$ of the form given in 
Lemma \ref{cc}.
Suppose first that $(X,Y,i,j)$
$\in \M'_k$, i.e. $k'=k''=k$.
In this case, it is easy to see that the  isotropy in
$G$ and $\g$ of the triple $(Y,i,j)$,
hence of the quadruple $(X,Y,i,j)$, are trivial. Moreover,
the choice of a representative of the form 
in Lemma \ref{cc}, in the same orbit of $(X,Y,i,j)$,
is unique up to the action of $S_k \times S_{n-k}$, where
$S_k$ permutes the first $k$ coordinates of $\C^n$ and
$S_{n-k}$ permutes the last $n-k$ coordinates of $\C^n$. 
Denote by $\Delta$ the big diagonal in $\C^n$. 
By Lemma \ref{cc}, we have
$$ \M'_k \simeq G\times_{S_{k}\times S_{n-k}} 
\big( \C^n \times (\C^n\backslash\Delta) \big),$$
where the action of $S_{k}$, resp. $S_{n-k}$, on 
$\C^n \times (\C^n\backslash\Delta) $ is by permutation of 
the pairs $(x_1,y_1)$, $\ldots$, $(x_k,y_k)$, resp.
$(x_{k+1},y_{k+1})$, $\ldots$, $(x_n,y_n)$.
This proves \vi.

Suppose now that $(X,Y,i,j)\in \M''_k$, so that 
$k''>k'$. In this case, the subgroup of $G$
consisting of matrices of the form 
$$\mathrm{diag}(\underbrace{1,\ldots, 1}_{k'}, a, b, \ldots, c,
\underbrace{1,\ldots, 1}_{n-k''})\quad \mbox{ where }
a,b,\ldots, c\in \C\backslash\{0\},$$ 
acts trivially on $(X,Y,i,j)$. This is a subgroup of strictly
positive dimension.
Hence, Lemma \ref{cc} implies \vii. 
\endproof

\begin{proof}[Proof of Theorem \ref{t1}.]
Let $\Delta$ be the big diagonal in $\C^{(n)}$.
By Lemma \ref{nak}, we have a set-theoretic equality
\begin{equation}\label{union}
\M=\left(\bigcup\limits_{k=0}^n \M'_k\right)
\sqcup \left(\bigcup\limits_{k=0}^n \M''_k\right)
\sqcup \pi^{-1}(\Delta).
\end{equation}

We claim that the decompositions above imply that 
$\M_0$, $\ldots$, $\M_n$ are precisely
all the irreducible components of $\M$.
To this end, observe that,
by Lemma \ref{dd}, $\M'_k$ is connected and has dimension
$n^2+2n$, and $\M''_k$ has dimension strictly less than
$n^2+2n$. Further,
since $\dim\Delta <n$, Corollary \ref{flat}
 yields
$\dim\pi^{-1}(\Delta)<n+(n^2+n)=n^2+2n.$

Now, by  Corollary \ref{flat}, the scheme $\M$ is a complete intersection.
Hence, each
irreducible component
of $\M$ must have dimension  $n^2+2n$.
Thus, parts (i) and (ii)  of the Theorem both
follow from  \eqref{union}. 

Finally, since the $\g$-action on $\M'_k$ is free, 
the moment map $\mu$ is a submersion at generic points of $\M$.
Hence $\M$ is generically reduced.
But $\M$ is a complete intersection, hence it is 
Cohen-Macaulay. It follows that $\M$ is reduced, cf. 
\cite[Theorem 2.2.11]{CG} or \cite[Exercise 18.9]{Ei}.
This proves (iii).
\end{proof}

\subsection{Around commuting variety.}\label{around}
Let $I_1\subset\C[\g\times\g]=\C[X,Y]$ be the ideal generated
by all the $2\times2$ minors of the matrix $[X,Y]$.
Recall also the ideals $J\subset \C[X,Y,i,j]$ and
$I\subset \C[X,Y]$ defined in \S\S\ref{int1}-\ref{int2}.
Further, let $p^*: \C[X,Y] \hookrightarrow \C[X,Y,i,j]$ be the 
pullback morphism induced by the projection map $p: (X,Y,i,j)\mapsto (X,Y)$.
It is clear from definitions that one has
$$ I\supset I_1\quad\text{and}\quad p^*(I_1)\subset J.
$$
Thus,  there are natural algebra morphisms 
$$\C[\zz]:=\C[X,Y]/I\longleftarrow
 \C[X,Y]/I_1\stackrel{p^*}\longrightarrow\C[X,Y,i,j]/J.
$$
By \cite[Theorem 12.1]{EG}, we know that $I^G=I_1^G$. Hence, taking
$G$-invariants, we obtain the following diagram
\begin{equation} \label{nn}
\C[\zz]^G=\bigl(\C[X,Y]/I\bigr)^G\stackrel{\sim}\longleftarrow
 \bigl(\C[X,Y]/I_1\bigr)^G\stackrel{p^*}\longrightarrow
\bigl(\C[X,Y,i,j]/J\bigr)^G.
\end{equation}

\begin{proposition} \label{p1}
The morphism $p^*$ in (\ref{nn}) is an isomorphism.
\end{proposition}
\proof
To prove surjectivity, we note by Weyl's fundamental theorem
of invariant theory that $\C[X,Y,i,j]^G$ is generated as an
algebra by polynomials of the form $\Tr(P(X,Y,ij))$, where $P$
is any noncommutative polynomial. But this is equal to
$p^*(\Tr(P(X,Y,-[X,Y])))$ modulo $J$, and 
$\Tr(P(X,Y,-[X,Y]))$ is contained in $\C[X,Y]^G$.

For injectivity, since we know that  $I^G=I_1^G$,
it suffices to show that any element
$f\in \C[X,Y]^G\cap J$ is contained in $I$.
From the definition of $J$, we can write $f=\sum_a f_ad_a$ 
where $f_a\in\C[X,Y,i,j]$ and $d_a=[X,Y]_{rs}+i_r j_s$
for some $r,s$. Let $\phi: \g\times\g\hookrightarrow
\g\times\g\times V\times V^*,$
$(X,Y)\mapsto (X,Y,0,0)$
be the natural imbedding, and
 $\phi^*:\C[X,Y,i,j]\too \C[X,Y]$ the restriction morphism,
so  $\phi^*(i)=\phi^*(j)=0$.
Since $f$ is contained in $\C[X,Y]$, we have 
$f= \phi^*(f) = \sum_a \phi^*(f_a)\phi^*(d_a)$. 
But $\phi^*(d_a)=[X,Y]_{rs}$.
Hence, $f$ is contained in $I$. 
\endproof

\proof[Proof of Theorem \ref{t2}.]
By Theorem \ref{t1}, $\M$ is reduced, that is, 
$J=\sqrt{J}$, and so $J^G=(\sqrt{J})^G$.
But $(\sqrt{J})^G$ is the radical of $J^G$
in $\C[X,Y,i,j]^G$.
By Proposition \ref{p1}, it follows that
$I^G$ is equal to its radical in $\C[X,Y]^G$,
hence $I^G = (\sqrt{I})^G$.
\endproof

\subsection{}\label{state}
We are going to state a corollary of Theorem \ref{t2}
that will be used later
in applications to Cherednik algebras.

To this end, let
$\h:=\C^n$ be the
 permutation representation of  the Symmetric group $S_n$,
and let $S_n$ act diagonally on
$\h\times\h$.
The quotient $(\h\times\h)/S_n$ has a natural
structure of algebraic variety, with
coordinate ring $\C[(\h\times\h)/S_n]=\C[\h\times\h]^{S_n}.$

We may identify $\h$ with the Cartan subalgebra in $\g$ formed by diagonal
matrices, so
we have a tautological imbedding $\h\times\h\into\g\times\g$.
Write $i_o$ for the vector $i_o:=(1,1,\ldots,1)\in V$.
We define the following closed imbedding
\begin{align}\label{j}\beps:\ &
 \h\times\h\into \g\times\g\times V\times V^*,\\
&(x_1,\ldots,x_n,\,y_1,\ldots,y_n)\mto
\big({\operatorname{diag}}(x_1,\ldots,x_n),\,
{\operatorname{diag}}(y_1,\ldots,y_n),\,i_o,\,0\big).
\nonumber
\end{align}
Note that
$S_n$, viewed as the subgroup of permutation matrices in $G$,
fixes $i_o$. Thus, the image of $\beps$ is an $S_n$-stable
subset in $\g\times\g\times V\times V^*$.

Let $\M_\red$ denote the scheme $\M$ with reduced scheme structure
(of course, we already know by Theorem \ref{t1}(iii) that $\M$ is
reduced, but we prefer not to use the Theorem at this point). 
It is clear that the above defined map $\beps$ gives
an $S_n$-equivariant closed imbedding $\beps: \h\times\h\into \M_\red$.

Further, it is a well known and easy
consequence of Weyl's fundamental theorem on $GL_n$-invariants
that restriction of polynomial functions from $\g\times\g$
to $\h\times\h$ induces an algebra isomorphism
$q: (\C[X,Y]/\sqrt{I})^G\iso \C[\h\times\h]^{S_n}.$
(Here and below we use the notation of \eqref{nn}).
Thus, using \eqref{nn}, we obtain
 the following chain of algebra isomorphisms
\begin{align}\label{composite}
\C[\M]^G=\bigl(\C[X,Y,i,j]/J\bigr)^G
&\cong\bigl(\C[X,Y]/I_1\bigr)^G\cong\bigl(\C[X,Y]/I\bigr)^G\\
&=
\C[X,Y]^G/I^G=
\C[X,Y]^G/(\sqrt{I})^G\stackrel{q}\cong\C[\h\times\h]^{S_n}.
\nonumber
\end{align}

Let $f: {\M_\red}\map (\h\times\h)/S_n$ be the morphism
of schemes induced by the composite  algebra isomorphism in 
\eqref{composite}. Set theoretically, the morphism
$f$ can be described in more geometric terms as follows.

Given an upper triangular matrix $X$,
let $X_{\operatorname{diag}}\in\h$ denote
the diagonal part of $X$.
Recall further that, for any quadruple
$(X,Y,i,j)\in\M_\red$, the matrices $X,Y$ can be
{\em simultaneously} put into upper triangular form,
by \cite[Lemma 12.7]{EG}. 
We assign to such a   quadruple
$(X,Y,i,j)\in\M_\red,$ where $X,Y$ are  upper triangular matrices,
the pair $(X_{\operatorname{diag}},Y_{\operatorname{diag}})\in\h\times\h$,
taken up to  $S_n$-diagonal action on $\h\times\h$.
The resulting map $\M_\red\map (\h\times\h)/S_n$ 
 is clearly
constant on $G$-diagonal orbits in $\M_\red.$
Furthermore, it is easy to verify that this
map of {\em sets} corresponds to the scheme morphism
$f$ defined earlier using the chain of  algebra 
isomorphisms in \eqref{composite}.

\begin{lemma}\label{f} Restriction of functions via $\beps$, resp.
pull-back of functions via $f$, induce mutually inverse  graded algebra
isomorphisms $\beps^*$ and $f^*$ in the diagram below:
$$
\xymatrix{
\C[\h\times\h]^{S_n}=\C[(\h\times\h)/S_n]
\ar@<1ex>[rr]^<>(0.5){f^*}_<>(0.5){\sim}&&
\C[\M_\red]^G\ar@<1ex>[ll]^<>(0.5){\beps^*}&&&\C[\M]^G.
\ar[lll]_<>(0.5){_{\text{Theorem \ref{t1}(iii)}}}^<>(0.5){\sim}
}
$$
\end{lemma}
\proof 
It is straightforward to verify
that $\beps^*\ccirc
f^*=\Id_{\C[\h\times\h]^{S_n}}$. Further, the map
$f^*$ is an isomorphism by construction, see
 \eqref{composite}.   
It follows that   $\beps^*$ is the inverse of $f^*$.
\endproof

\section{Generalization to quiver moment maps}
\subsection{Quiver setting.}
Throughout this section, we let $Q$ be a quiver with vertex set $I$.
The double $\QQ$ of $Q$ is the quiver obtained from $Q$ by
adding a reverse edge $a^*: {j\to i}$ 
for each edge $a:{i\to j}$ in $Q$.
If $a: {i\to j}$ is an edge
in $\QQ$, we call $t(a):=i$ its tail, and 
$h(a):=j$ its head.
The opposite quiver $Q^{op}$ is the quiver with the same
underlying graph as $Q$ but with all the edges oriented
in the opposite direction to the ones in $Q$.

On $\C^I$, we have the standard inner product
$\al\cdot\beta:=\sum_{i\in I} \al_i\beta_i,$ and
 we write  $|\al|^2:=\sum_{i\in I} \al_i^2.$
We will also use
the {\it Ringel form} of $Q$, a (not necessarily symmetric) bilinear form on $\Z^I$
defined by
$$\langle \al,\beta \rangle
:= \sum_{i\in I} \al_i\beta_i
- \sum_{a\in Q} \al_{t(a)}\beta_{h(a)}, \quad\mbox{ where }
\alpha=(\al_i)_{i\in I},
\ \beta=(\beta_i)_{i\in I}.$$
Let $(\al,\beta):= \langle \al,\beta \rangle
+\langle \beta,\al \rangle$ be its symmetrization.
The corresponding quadratic form $q(\al):=
\langle \al,\al \rangle = \frac{1}{2}(\al,\al)$
is the {\it Tits form}; we set
$$p(\al) := 1-q(\al)
= 1 + \sum_{a\in Q} \al_{t(a)}\al_{h(a)} - |\al|^2.$$

Let $\epsilon_i\in \Z^I$ denote the coordinate vector
corresponding to the vertex $i\in I$.

The representations of $Q$ of dimension
vector $\al\in \N^I$ are the elements of the 
vector space 
$$ \Rep(Q,\al) := \bigoplus_{a\in Q} 
\Mat(\al_{h(a)}\times \al_{t(a)}, \C). $$
Let 
$G(\al):= \big(\prod_{i\in I} GL(\al_i, \C)\big)/\C^*$. 
The group $G(\al)$
acts on $\Rep(Q,\al)$ by conjugation and the orbits are
the isomorphism classes of representations of $Q$ of
dimension vector $\al$.
Let $$\End(\al)_0 := \big\{ (g_i)_{i\in I} \mid 
\sum_{i\in I} \Tr (g_i) =0 \big\}
\subseteq \End(\al):=\bigoplus_{i\in I} \Mat(\al_i, \C).$$

We identify $\Rep(\QQ,\al)$ 
with the cotangent bundle
$T^* \Rep(Q^{op},\al)$, and $\End(\al)_0$ with the dual of the Lie
algebra of $G(\al)$. Then, we have a moment map
$$ \mu_\al: \Rep(\QQ,\al) \too \End(\al)_0,
\qquad \mu_\al(x) := \sum_{a\in Q} [x_a, x_{a^*}].$$
Let $\la\in \C^I$. The fiber $\mu_\al^{-1}(\la)$ is 
a (not necessarily reduced) scheme.
In particular, the zero fiber, $\mu_\al^{-1}(0),$ is the
union (possibly infinite) 
 of the conormal bundles to $G(\alpha)$-orbits in  $\Rep(Q^{op},\al)$.

\subsection{Irreducible components of $\mu_\al\inv(0)$.}
 Fix $\la\in \C^I$, and let $R_\la^+$ be the set of positive roots
$\al$ with $\al\cdot \la=0$.
Let $\Sigma'_\la,$ resp., $\Sigma_\la$, be the set of $\al\in R_\la^+$ with the
property that 
\begin{align}\label{Sigma}
&p(\al) \geq  p(\be^{(1)}) +\cdots+ p(\be^{(r)}),\quad
\text{resp.},
\quad p(\al) > p(\be^{(1)}) +\cdots+ p(\be^{(r)}),\\
&\text{for any decomposition}\quad\al = \be^{(1)}+\cdots+\be^{(r)}\quad
\text{where}\quad\be^{(t)}\in R_\la^+, \,\forall t=1,\ldots,r,
\quad\text{and}\quad r\geq 2.\nonumber
\end{align}

Given $\al\in \Sigma'_\la$,  let
$\Sigma'_\la(\al)$ be the set of 
decompositions $\al = \be^{(1)}+\ldots+\be^{(r)},\,\be^{(t)}\in R_\la^+$
such that the inequality in \eqref{Sigma} is an equality,
i.e., such that we have
 $p(\al) = p(\be^{(1)}) +\cdots+ p(\be^{(r)})$.

\begin{theorem}\label{irr_thm}
Assume $\al \in \Sigma'_\la$. Then 

\vi The scheme
$\mu_\al^{-1}(\la)$  is equidimensional of dimension
$|\al|^2 -1 +2p(\al)$; furthermore, it is a
complete intersection in $\Rep(\QQ,\al)$.

\vii 
The irreducible components of $\mu_\al^{-1}(\la)$ 
are in 1-1 correspondence with elements of the set $\Sigma'_\la(\al)$. 
\end{theorem}

Note that  $\Sigma_\la \subset \Sigma'_\la$. Also, if $\al \in \Sigma_\la$, then
$\Sigma'_\la(\al) = \{ \al \}$. 
Crawley-Boevey
proved
in \cite[Theorem~1.2]{CB} that, for any  $\al\in \Sigma_\la$,
the scheme  $\mu_\al^{-1}(\la)$ is reduced and {\em irreducible}.

\begin{conjecture}\label{red_conj} For any $\al\in \Sigma'_\la$, the scheme
 $\mu_\al^{-1}(\la)$ is  reduced.
\end{conjecture}

\begin{proof}[Proof of theorem \ref{irr_thm}] Part (i) of the Theorem
is   \cite[Theorem 4.4]{CB}.

To prove (ii) we need to introduce some notation.
Given an algebraic group $G$ acting on an algebraic variety
$Z$, let $\dim_GZ$ denote the maximal number 
of parameters for $G$-orbits in $Z$, see \cite[\S3]{CB}.

 For any integer $d\geq 1$ and dimension vector $\be$,
let $I^d(\be)$ be the set of indecomposable
 representations $\rho\in\Rep(Q^{op},\be)$,
 such that $\dim G(\be)\cdot\rho=d$. It is known, see 
 \cite[pp.142-143]{KR}, that
 each set $I^d(\be)$ has the structure of a locally closed reduced
subscheme in $\Rep(Q^{op},\be).$  Thus,
$I(\be)=\sqcup_{d\geq 0} I^d(\be)$ is the constructible
set of  indecomposable
 representations of $Q^{op}$ of dimension $\be$.

Given an $r$-tuple of dimension vectors
$\be^{(1)}, \ldots, \be^{(r)}\in\Z^I$, 
we consider the action of the group
$H:=G(\be^{(1)})\times\ldots\times G(\be^{(r)})$ on
the set $
I(\be^{(1)})\times\ldots\times
I(\be^{(r)}).$ It is known that
$\dim_H\bigl(I(\be^{(1)})\times\ldots\times
I(\be^{(r)})\bigr)\leq \sum_{t=1}^r p(\be^{(t)}),$
by 
Kac's theorem, \cite{Ka}. Furthermore,
it was pointed out  to us by
Crawley-Boevey that,  according to the last remark
in \cite[p.144]{KR}, one has

\begin{claim}\label{bk} {\em For any  $r$-tuple $\be^{(1)}, \ldots,
\be^{(r)}$
of positive roots,
there exists exactly one $r$-tuple   $(d_1,\ldots,d_r)$ and exactly one
 irreducible component,
$Z(\be^{(1)}, \ldots, \be^{(r)})$,
of the set $I^{d_1}(\be^{(1)})\times\ldots\times
I^{d_r}(\be^{(r)})$ such that $\dim_HZ(\be^{(1)}, \ldots, \be^{(r)})
=\sum_{t=1}^r p(\be^{(t)}).$ 

Moreover, for any $H$-stable
constructible set $\mathcal{Y}\sset I(\be^{(1)})\times\ldots\times
I(\be^{(r)})$, we have}
\beq{strict}
\dim_H\mathcal{Y} <
\dim_HZ(\be^{(1)}, \ldots, \be^{(r)})\en\;
\text{whenever}\en\; \dim\bigl(\mathcal{Y}
\cap Z(\be^{(1)}, \ldots, \be^{(r)})\bigr)<\dim
Z(\be^{(1)}, \ldots, \be^{(r)}).
\eeq
\end{claim}

Now, fix a dimension vector $\al$ and a
 decomposition $\al=\be^{(1)}+\cdots+\be^{(r)}.$
We have
 a  natural (block diagonal) imbedding
$H:=G(\be^{(1)})\times\ldots\times G(\be^{(r)})\into
G(\al),$ and a similar imbedding
$$
I(\be^{(1)})\times\ldots\times
I(\be^{(r)})\sset \Rep(Q^{op},\be^{(1)})\times\ldots\times
\Rep(Q^{op},\be^{(r)})\into \Rep(Q^{op},\al).
$$
Write $I(\be^{(1)}, \ldots, \be^{(r)})\sset  \Rep(Q^{op},\al)$
for the set of representations whose indecomposable summands
have dimensions $\be^{(t)},\,t=1,\ldots,r.$
Equivalently,
we have  $$I(\be^{(1)}, \ldots, \be^{(r)}):=
G(\al)\cd\bigl(I(\be^{(1)})\times\ldots\times
I(\be^{(r)})\bigr),$$
is the
$G(\al)$-saturation of the set $I(\be^{(1)})\times\ldots\times
I(\be^{(r)}).$
It is clear that the space $\Rep(Q^{op}, \al)$ is a disjoint
union
of the sets $I(\be^{(1)}, \ldots, \be^{(r)})$ for various decompositions
of $\al=\be^{(1)}+\cdots+\be^{(r)}$ such that  each $\be^{(t)}\in
R_\la^+$.

Let $\pi:\mu_\al^{-1}(\la)\into\Rep(\QQ, \al)\to\Rep(Q^{op}, \al)$ be the 
restriction of the vector bundle projection
$T^*\Rep(Q^{op}, \al)$
$\to\Rep(Q^{op}, \al).$
Thus, the set $\mu_\al^{-1}(\la)$ breaks up into a disjoint
union
of various pieces $\pi^{-1}\bigl(I(\be^{(1)}, \ldots,
\be^{(r)})\bigr)$.
Fix one such piece corresponding to
 a  decomposition $\al = \be^{(1)}+\cdots+\be^{(r)},$
where each $\be^{(t)}\in R_\la^+$.
It follows from 
 \cite[Lemma 3.4]{CB} and \cite[Lemma 4.3]{CB} that
 $\dim \pi^{-1}\bigl(I(\be^{(1)}, \ldots,
\be^{(r)})\bigr)
\leq |\al|^2 -1 + 2p(\al)$ with equality
if and only if
the decomposition $(\be^{(1)},\ldots,\be^{(r)})$
is in $\Sigma'_\la(\al)$. 
Therefore, by part (i) of the theorem, we see that  any irreducible component
of $\mu_\al^{-1}(\la)$ has to be the closure of an
irreducible component of the set  $\pi^{-1}\bigl(I(\be^{(1)}, \ldots,
\be^{(r)})\bigr)$ such that
 the  decomposition $\al = \be^{(1)}+\cdots+\be^{(r)}$
belongs to  $\Sigma'_\la(\al)$; moreover, this irreducible
component 
must have dimension equal to $|\al|^2 -1 + 2p(\al).$

To complete the proof of part (ii) of the Theorem,
it remains to show that, for each $r$-tuple
$(\be^{(1)},\ldots,\be^{(r)})\in\Sigma'_\la(\al)$, the
set $\pi^{-1}\bigl(I(\be^{(1)}, \ldots,
\be^{(r)})\bigr)$ contains only one irreducible component
of required dimension. To this end, 
recall the  set
$Z(\be^{(1)}, \ldots, \be^{(r)})$ introduced earlier.
This is an $H$-stable irreducible subvariety
in $\Rep(Q^{op},\al)$. Let 
$\bd$
 denote the dimension of the $G(\al)$-orbit through a general point of
 $Z(\be^{(1)}, \ldots, \be^{(r)})$.
The set $Z^{\text{reg}}(\be^{(1)}, \ldots, \be^{(r)})$
of all points  $\rho\in Z(\be^{(1)}, \ldots, \be^{(r)})$
such that $\dim G(\al)\cdot\rho=\bd$ is a Zariski
open dense subset. Using  \cite[Lemma 3.4]{CB} and \cite[Lemma 4.1]{CB}
we find that
 $\dim \pi^{-1}\bigl(G(\al)\cdot Z^{\text{reg}}(\be^{(1)}, \ldots,
\be^{(r)})\bigr)= |\al|^2 -1 + 2p(\al)$.
Furthermore, since all  $G(\al)$-orbits in 
$G(\al)\cdot Z^{\text{reg}}(\be^{(1)}, \ldots,
\be^{(r)})$ have the same dimension, we deduce that
all fibers of the projection
$\pi^{-1}\bigl(G(\al)\cdot Z^{\text{reg}}(\be^{(1)}, \ldots,
\be^{(r)})\bigr)\map
G(\al)\cdot Z^{\text{reg}}(\be^{(1)}, \ldots,
\be^{(r)})$ are affine-linear spaces of the same dimension,
hence, the set $\pi^{-1}\bigl(G(\al)\cdot Z^{\text{reg}}(\be^{(1)}, \ldots,
\be^{(r)})\bigr)$ is irreducible.

On the other hand, put $\mathcal{Y}:=\bigl(I(\be^{(1)})\times\ldots\times
I(\be^{(r)})\bigr)
\sminus Z^{\text{reg}}(\be^{(1)}, \ldots, \be^{(r)}).$
Then, from \eqref{strict}  and \cite[Lemma 4.1]{CB},
we obtain
 $$\dim_{G(\al)}G(\al)\cd\mathcal{Y}=\dim_H\mathcal{Y}<
\dim_HZ(\be^{(1)}, \ldots, \be^{(r)})=\sum\nolimits_{t=1}^r
p(\be^{(t)}).$$

Thus, using 
 \cite[Lemma 3.4]{CB} we deduce that
$$\dim \pi^{-1}\bigl(I(\be^{(1)}, \ldots,
\be^{(r)})\sminus
G(\al)\cdot Z^{\text{reg}}(\be^{(1)}, \ldots,
\be^{(r)})\bigr)=
\dim \pi^{-1}\bigl(G(\al)\cdot\mathcal{Y}\bigr)< |\al|^2 -1 + 2p(\al),
$$
and the theorem follows.
\end{proof}

\begin{remark} We do not know if 
$\pi^{-1}(I(\be^{(1)}, \ldots, \be^{(r)}))$
is irreducible, for any 
$(\be^{(1)}, \ldots, \be^{(r)})\in \Sigma'_\la(\al)$.
%
\end{remark}

\subsection{Extended Dynkin case.}
One of the goals of this section is to prove Theorem \ref{tquiver}.
Using the McKay correspondence, one may reformulate
the theorem in the language of quivers, see eg. \cite[\S11]{EG}.
This will allow us to use the results from \cite{CB}.

From now on, for the rest of this section, 
we let $Q$ be an affine Dynkin quvier,
let $o$ be an extending vertex of $Q$, and let
$S$ be the quiver obtained from $Q$ by adjoining one 
vertex $s$ and one arrow $s\to o$. 

Fix a positive integer $n$.
Let $I$ be the vertex set of $S$, let $\delta\in \N^I$ be
the minimal positive imaginary root of $Q$, 
and let $\al:=n\delta+\epsilon_s$.
We have $q(\delta)=0$ and $q(\al)=1-n$,
so $p(\delta)=1$ and $p(\alpha)=n$.

We will use some of the notation introduced in the proof of Theorem
\ref{irr_thm}. In particular,
for any decomposition $\al=\be^{(1)} +\cdots+ \be^{(r)},$ 
where $\be^{(t)}\in \N^I\setminus\{0\}$,
let $I(\be^{(1)},\ldots,\be^{(r)})$ be the subset of 
$\Rep(S^{op},\al)$ consisting of the representations of
$S^{op}$ whose indecomposable summands have dimension $\be^{(t)},\,t=1,\ldots,r.$

Fix  $\la\in \C^I$ such that $\la_s = 0$ and $\la\cdot\delta=0$.
We have the moment map
$$\mu_\al: \Rep(\SS,\al) \too \End(\al)_0.$$
Denote by $\pi:\mu_\al^{-1}(\la) \too \Rep(S^{op},\al)$
the projection map.

\begin{lemma}  \label{cb1}
For any decomposition $\al=\be^{(1)} +\cdots+ \be^{(r)}$
with $\be^{(t)}\in \N^I\setminus\{0\}$, we have
$$\dim \pi^{-1}(I(\be^{(1)},\ldots,\be^{(r)}))
\leq n^2|\delta|^2+2n$$ with equality if and only if
all but one of the $\be^{(t)}$ are equal to $\delta$.
\end{lemma}
\begin{proof}
The lemma follows from \cite[Lemma 9.2]{CB},
\cite[Lemma 4.3]{CB}, and \cite[Lemma 3.4]{CB}.
\end{proof}

Let $\sigma: \Rep(S^{op},\al) \too \Rep(Q^{op}, n\delta)$ be
the projection map.
Let $U \subset \Rep(Q^{op}, n\delta)$ be the subset consisting 
of all representations $Y$ which have a decomposition 
$Y\cong Y_1\oplus\cdots\oplus Y_n$, where each $Y_t$ is indecomposable
with dimension vector $\delta$, and $\dim \End(Y) =n$
(this implies in particular that  $\End(Y_t)=\C,$ for all $t=1,\ldots,n$).
It is well known that the canonical decomposition of
$n\delta$ is $\delta+\cdots+\delta$, 
and moreover, the subset
$\{Y \mid \dim \End(Y)\leq n'\}$ is open 
in $\Rep(Q^{op}, n\delta)$,
for any $n'$; see \cite{Ka} and \cite{Sc}.
Hence, $U$ is a dense open subset of $\Rep(Q^{op}, n\delta)$.

Now, we have the following composition of maps:
$$
\Rep(\SS,\al)
\stackrel{\pi}{\too} \Rep(S^{op}, \al)
\stackrel{\sigma}{\too} \Rep(Q^{op}, n\delta).$$

\begin{lemma} \label{cb2}
We have $\dim \pi^{-1}\big( \sigma^{-1} ( 
\Rep(Q^{op}, n\delta) \setminus U ) \big) <
n^2|\delta|^2 + 2n$.
\end{lemma}
\begin{proof}
This is immediate from \cite[Lemma 9.3]{CB}
and \cite[Lemma 3.4]{CB}.
\end{proof}

A representation of a quiver is called a brick if 
its endomorphism algebra is of dimension one.
Let $\M'_k$ (resp., $\M''_k$)
be the subset of $\pi^{-1}( I(k\delta+\epsilon_s,
\delta,\ldots,\delta) )$ consisting of all $$(X,Y,i,j)\in
\Rep(Q,n\delta)\times\Rep(Q^{op},n\delta)\times V\times V^*$$
such that $Y\in U$, and $(X,Y,i,j)$ is a brick
(resp., $(X,Y,i,j)$ is not a brick).
Define $\M_k$ to be the closure of $\M'_k$ in 
$\mu_\al^{-1}(\la)$.

The main result of this section is Theorem \ref{quiverversion}
below, which is a
 more precise version of Theorem \ref{tquiver}.
Note  that, for the quiver $S$ and $\al$, $\la$ defined above,
in the notation of \eqref{Sigma} we have $\al\in \Sigma'_\la,$ 
by \cite[Lemma 9.2]{CB}. Thus, part (i) in  the theorem below
is a special case of Theorem \ref{irr_thm}(i) while part (ii) 
in  the theorem below provides an explicit description of the
1-1 correspondence from Theorem \ref{irr_thm}(ii).

Our proof of Theorem \ref{quiverversion} will be independent
of the proof of  Theorem \ref{irr_thm}. We remark also that,
most important in  the theorem below,  is its part (iii),
which is a special case of  Conjecture \ref{red_conj}.

\begin{theorem} \label{quiverversion}
\vi $\mu_\al^{-1}(\la)$ is  equidimensional,
we have $\dim \mu_\al^{-1}(\la) = n^2|\delta|^2+2n$;
furthermore, the scheme $\mu_\al^{-1}(\la)$ is a complete intersection
in $\Rep(\SS,\al)$.

\vii The irreducible components of $\mu_\al^{-1}(\la)$ are
$\M_0, \ldots, \M_n$.

\viii The scheme $\mu_\al^{-1}(\la)$ is reduced.
\end{theorem}
First, recall the following standard result, cf. \cite[Lemma 3.1]{CB}.

\begin{lemma} \label{lifting}
If $y=(y_{a^*})_{a\in Q} \in \Rep(Q^{op},\al)$,
then there is an exact sequence
$$ 
0 \too \Ext^1(y,y)^* \stackrel{\mathsf r}\too \Rep(Q,\al) 
\stackrel{\mathsf c}{\too} \End(\al) 
\stackrel{\mathsf t}{\too} \End(y)^* \too 0,
$$
where the map $\mathsf c$ sends $(x_a)\in \Rep(Q,\al)$
to $\sum_{a\in Q} [x_a, y_{a^*}]$,  the map
$\mathsf t$ sends $(g_i)$ to the linear map
$\End(y)\to \C: (\phi_i)\mapsto \sum_i \Tr(g_i\phi_i),$
and
the map $\mathsf r$ was defined in
\cite[\S3]{CB}. \qed
\end{lemma}

Next, we prove the following lemma.

\begin{lemma} \label{cb3}
\vi We have $\dim \M'_k = n^2|\delta|^2 + 2n$, 
and $\M'_k$ is irreducible.

\vii We have $\dim \M''_k < n^2|\delta|^2 + 2n$.
\end{lemma}
\begin{proof}
Let $R=\Rep(Q^{op},\delta)\times 
\cdots \times \Rep(Q^{op},\delta)$,
and consider it as a subset of $\Rep(Q^{op},n\delta)$ using
block-diagonal matrices.
Let $I'_k$ be the subset of $\Rep(S^{op},\al)$ 
consisting of the elements $(Y,j)$ such that
$Y\in R\cap U$ and
$j=(\underbrace{1,\ldots,1}_{k}, 0, \ldots, 0)$.
We have $I'_k\subset I(k\delta+\epsilon_s,\delta,\ldots,\delta)$.
It is clear that
$\pi^{-1}(I'_k)$ is contained in the disjoint union 
$\M'_k \sqcup \M''_k$. 
Moreover, any element in $\M'_k \sqcup \M''_k$
is $G(\al)$-conjugate to an element in $\pi^{-1}(I'_k)$.
Therefore, the $G(\al)$-saturation of $\pi^{-1}(I'_k)$ is
$\M'_k \sqcup \M''_k$.

Now let $y=(Y,j)\in I'_k$, and write
$Y=Y_1\oplus\cdots\oplus Y_n$,
where $Y_t\in \Rep(Q^{op},\delta)$.
Observe that $\dim \End(y) = n-k+1$, hence by
Lemma \ref{lifting}, the fiber
$\pi^{-1}(y)$ is an affine space of dimension
$2n-k$. 
Since $I'_k$ is an irreducible,
it follows that the scheme $\pi^{-1}(I'_k)$ is 
irreducible and has dimension $n|\delta|^2 + 2n-k$.
Hence, $\M'_k \sqcup \M''_k$ is irreducible.

We begin the proof of the inequality of part (ii) of the lemma.
Let $p: \pi^{-1}(y) \too V$ be the 
projection map $p(X,Y,i,j)=i$.

\medskip

{\it Claim:} The vector $i$ is in the image of $p$ 
if and only if it is of the form
$i=(0,\ldots,0, i_{k+1},\ldots, i_n)$.
If $i$ is in the image of $p$, then $p^{-1}(i)$ is 
an affine space of dimension $n$.

\medskip

{\it Proof of Claim:}
Let $N=\sum_t \delta_t$.
We shall consider $X$ and $Y$ as endomorphisms of a 
vector space of dimension $nN$.
Recall that we write $Y$ as a block-diagonal matrix 
$Y_1\oplus\cdots\oplus Y_n$. Similarly, we
write $X\in \Rep(Q,n\delta)$ 
as a $n\times n$ block matrix, where each
block is a $N\times N$ matrix.

Now consider Lemma \ref{lifting} for the quiver $Q$
and dimension vector $n\delta$.
Since $Y\in U$, we have $\dim\End(Y)=n$, so
the cokernel of the map $\mathsf c$ in Lemma \ref{lifting}
has dimension $n$.
For any $X\in \Rep(Q,n\delta)$, the element
$\mathsf c(X)=[X,Y]$ is a $n\times n$ block matrix
such that each of the $n$ blocks on the diagonal 
have trace $0$.
Hence, the image of $\mathsf c$ is the subspace of
all $n\times n$ block matrix
such that each of the $n$ blocks on the diagonal
have trace $0$.
Since $\la\cdot\delta=0$ and 
$j=(\underbrace{1,\ldots,1}_{k},0,\ldots,0)$,
the element $\la-ij$ is in the image of $\mathsf c$ if and only if
$i$ is of the form $i=(0,\ldots,0,i_{k+1},\ldots, i_n)$,
and in this case, by Lemma \ref{lifting},
$p^{-1}(i)$ is an affine space of dimension $n$.

This completes the proof of the Claim.

\medskip

We can now prove that $\dim \M''_k < n^2|\delta|^2 + 2n$.
Suppose $(X,Y,i,j)\in \pi^{-1}(I'_k)$.
By the above Claim, we have 
$i=(0,\ldots, 0, i_{k+1}, \ldots, i_n)$.
If $i_{k+1},\ldots, i_n$ are all nonzero, then
it is clear that $(X,Y,i,j)$ is a brick.
Define $Z_k$ to be the closed subset of $\pi^{-1}(I'_k)$
defined by the equation 
$i_{k+1}i_{k+2}\cdots i_n=0$
(that is, one of the last $n-k$ entries of $i$ is
equal to $0$).
Then $\M''_k$ is contained in the $G(\al)$-saturation 
of $Z_k$.
By the above Claim, $Z_k$ is properly 
contained in $\pi^{-1}(I'_k)$,
so we have $\dim Z_k < n|\delta|^2 +2n-k$.

Let $GL(n\delta):= \prod_t GL(n\delta_t)$, and
let $H_k\subset GL(\delta)^n$ 
be the subgroup of $GL(n\delta)$ consisting of
diagonal block matrices whose component at the 
extending vertex $o$ is of the form $\mathrm{diag}
(1,\ldots, 1, g_{k+1}, \ldots, g_n)$.
We have $\dim H_k = n|\delta|^2 -k$.
Moreover, $Z_k$ is stable under the action $H_k$.
Hence,
\begin{align*}
\dim \M''_k \leq & \dim Z_k + \dim GL(n\delta) -\dim H_k \\
< & (n|\delta|^2 +2n-k) + n^2|\delta|^2
- (n|\delta|^2 -k) =  n^2|\delta|^2 +2n.
\end{align*}
Therefore, by Lemma \ref{cb1} and Lemma \ref{cb2},
it follows that $\M'_k$ is 
of dimension $n^2|\delta|^2 +2n$.
Moreover, $\M'_k$ is open dense in 
the $G(\al)$-saturation of $\pi^{-1}(I'_k)$,
hence it is irreducible.
\end{proof}

{\bf Proof of Theorem \ref{quiverversion}}:
Observe that $\mu_\al^{-1}(\la)$ is defined by $n^2|\delta|^2$
equations in the $2n^2|\delta|^2 +2n$ dimensional
space $\Rep(\SS,\al)$. Thus, the irreducible components
of $\mu_\al^{-1}(\la)$ must have dimension of at least 
$n^2|\delta|^2 +2n$.
By Lemma \ref{cb1}, it follows that $\mu_\al^{-1}(\la)$
is a complete intersection and equidimensional of dimension
$n^2|\delta|^2 +2n$.
Further, set theoretically, we have
$$ \mu_\al^{-1}(\la) =
\left(\bigcup_{k=0}^n \M'_k\right) \sqcup 
\left(\bigcup_{k=0}^n \M''_k\right)\sqcup
\pi^{-1}\big( \sigma^{-1}( \Rep(Q^{op},n\delta)\setminus U)\big).
$$
It follows, by Lemma \ref{cb2} and Lemma \ref{cb3},
that the irreducible components of $\mu_\al^{-1}(\la)$
are $\M_0, \ldots, \M_n$.

Since the representations in $\M'_k$ are bricks,
the moment map $\mu_\al$ is a submersion at 
generic points of $\mu_\al^{-1}(\la)$, and so
$\mu_\al^{-1}(\la)$ is generically reduced.
Moreover, since $\mu_\al^{-1}(\la)$ is a complete
intersection, it is Cohen-Macaulay.
Hence, $\mu_\al^{-1}(\la)$ is reduced 
(cf. \cite[Theorem 2.2.11]{CG} or \cite[Exercise 18.9]{Ei}).  \qed

\section{A Lagrangian variety}
\subsection{}
Let $\P:=\P(V)$ be the projective space (of dimension
$n-1$). We identify the total space of the cotangent bundle to $\P$
with 
\begin{equation}\label{TP}
T^*\P=\{(i,j)\in (V\sminus \{0\})\times V^*
\mid \langle j,i\rangle=0\}\big/\C^\times,
\end{equation}
where the multiplicative group $\C^\times$ acts
naturally on $V\times V^*$ by $t(i,j)= (t\cdot i,t\inv\cdot j)$.
The group $G=\GL(V)$ acts naturally on $\P$
and the induced $G$-action on $T^*\P$ is Hamiltonian
with moment map $(i,j)\mto ij\in\g=\g^*$.

Next, we set $\X:=\g\times\P$ and view it as a $G$-variety
with respect to  $G$-diagonal action. We remark that this
action  clearly
factors through the quotient $PGL(V)=GL(V)/\C^\times$,
in particular, any $G$-orbit in $\X$ may also be
regarded as either $SL(V)$- or $PGL(V)$-orbit.
Further, we have
\begin{equation}\label{TX}
T^*\X=T^*\g\times T^*\P=\{(X,Y,i,j)\in \g\times\g\times (V\sminus \{0\})\times V^*
\mid \langle j,i\rangle=0\}\big/\C^\times,
\end{equation}
where $\C^\times$ acts  on $(i,j)$ as above and does not act on $X,Y$.
 Again,
the induced $G$-action on $T^*\X$  is Hamiltonian
with  moment map given, essentially, by formula \eqref{moment}.

\subsection{A stratification of $\X$.}
 Given a  direct sum decomposition  $V=V_1\oplus\ldots\oplus V_l,$
write $L=\GL(V_1)\times\ldots\times\GL(V_l)$
for the corresponding Levi subgroup in $G=\GL(V),$ formed by
 block-diagonal matrices,
Let
$\fl=\Lie(L) =\gl(V_1)\oplus\ldots\oplus \gl(V_l)$
be the corresponding Levi subalgebra in $\g=\gl_n$
and $\mn_\fl=\mn\cap\fl$  the
nilpotent variety of the reductive Lie algebra $\fl$.
The group $L$ acts on $\mn_\fl$ with finitely many orbits.

We have a direct sum decomposition
$\fl=\fz_\fl\oplus[\fl,\fl],$
where 
$\fz_\fl$, the center of the Lie
algebra $\fl$, may be identified naturally with $\C^l$.
Let $\fzr_\fl\sset\fz_\fl=\C^l$
denote  the Zariski-open dense subset
formed by the elements of $\C^l$
with pairwise distinct coordinates.

G. Lusztig introduced, see
\cite{Lu1}, a certain stratification 
of the Lie algebra of  an
arbitrary reductive group $G$.
The strata of  Lusztig's stratification are smooth locally-closed
 $\Ad G$-stable subvarieties
labelled
by 
the set of $G$-conjugacy classes of  pairs $(\fl,\OO),$ where 
$\fl\sset\Lie G$ is a Levi subalgebra and $\OO\sset\fl$
is a {\em nilpotent}  $\Ad L$-conjugacy  class.

We are going to introduce a similar stratification
of the variety $\X$. 
The strata of our  stratification 
will be labelled by the set  $\scal$ formed $G$-conjugacy classes
of pairs
$(\fl,\Om)$,
where $\fl\sset\g$ is a Levi subalgebra
(with the corresponding Levi subgroup $L\sset G$),
and $\Om$ is an $L$-diagonal orbit
in  $\mn_\fl\times\P$.
Given such a pair $(\fl,\Om)\in\scal$,
we define
$$
\X(\fl,\Om):=\big\{(X, \ell)\in \g\times\P\mid(X,\ell)\in G(z+x,\C
v),\enspace \text{for
some}
\enspace z\in\fzr_\fl,\,
(x,v)\in\Om\big\},
$$
(here and below, we  write $G(X,\ell)$ for the $G$-diagonal orbit of
an element $(X,\ell)\in\g\times\P$;
we also use similar notation for  $L$-diagonal orbits).

One can show\footnote{This has been done by Rupert Yu
by generalizing Broer results
\cite{Bro} on the classical Jordan classes.}
that $\X(\fl,\Om)$ is a  smooth
locally-closed  $G$-stable subvariety of $\X$. Observe also that
each set  $\X(\fl,\Om)$ is stable under the
$\C^\times$-action on $\X=\g\times\P$ by dilations
along the $\g$-factor.

\begin{proposition}\label{finite} \vi  Any two pieces 
$\X(\fl,\Om)$ and $\X(\fl',\Om')$ are either equal or disjoint.

\vii The set $\scal$ is finite;  we have a finite decomposition
$\X=\bigsqcup_{(\fl,\Om)\in\scal}\,\X(\fl,\Om)$.
\end{proposition}

 \proof                                                                          
To prove (i), let $V=\oplus_r\,V_r$ and $V=\oplus_s\,V'_s$
be two direct sum decompositions, let $L$ and $L'$
be the corresponding Levi subgroups with Lie algebras
$\fl=\oplus_r\, \gl(V_r)$ and $\fl'=
\oplus_s\, \gl(V'_s)$, respectively.
Let $z\in \fzr_\fl,
z'\in\fzr_{\fl'},$ and let $(x,v)\in \mn_\fl\times V,\,
(x',v')\in \mn_{\fl'}\times V.$ 
Write $\Om=L(x,v),\, \Om'=L'(x',v')$,
for the corresponding orbits. 

Suppose there exists  $g\in G$ such that
 $(z'+x',\C  v')=g(z+x,\C  v)$.
By uniqueness of Jordan decomposition, we must have
$gzg^{-1}\in \fzr_{\fl'}$ and $gxg^{-1}\in\Ad L'(x')$.
Thus, the decompositions $V=\oplus_r\,V_r$ and $V=\oplus_s\,V'_s$
have the same number of direct summands,
moreover, we have
$g=wa$, where $w$ is a permutation matrix such that
each block of $\fl$ is mapped to a block of $\fl'$,
and $a\in L$.
Hence, the action of $g$ must take the set
$\{(z_1+x_1)\mid z_1\in \fzr_\fl,\, (x_1,\C  v_1)\in\Om\}$
into the set $\{(z_2+x_2)\mid z_2\in \fzr_{\fl'},\, (x_2,\C  v_2)\in\Om'\}$.
Part (i) follows.
                                  
To prove \vii, we observe that 
 there are only finitely many $\Ad L$-orbits
$\OO\sset\mn_\fl$. But for such an orbit $\OO$,
 the number
of $ L$-diagonal orbits in $\OO\times\P$
is {\em finite}, due to Lemma \ref{ee}.
Further, the centralizer in $\GL(V)$ of any $z\in \fzr_\fl$ is
$ L$. Part (ii) follows, since
it is clear by Jordan normal form
that any element of $\X$ belongs to some $\X(\fl,\Om)$.
 \endproof

Thus, $\dis\bigsqcup\nolimits_{(\fl,\Om)\in\scal}\, T^*_{\X(\fl,\Om)}\X$
is a (reducible) singular Lagrangian subvariety in $T^*\X$.

\subsection{Relevant strata.}\label{regstrata} We remind the reader  that the following
properties of
a linear map $X: V\to V$ are equivalent:

\pb{ The map $X$
has a cyclic vector, i.e., there is a line $\ell\in\P$ such that
$\C[X]\ell=V$;}

\pb{ $X$ is a
{\em regular} element of $\g$, i.e., such that $\dim \g^X=n$;}

\pb{ In Jordan normal form for $X$, different Jordan blocks
have pairwise distinct diagonal entries.}
 \smallskip

\begin{definition}\label{distinguished}
A pair $(X, \ell)\in\g\times\P$ is said to be 
{\em  relevant} if $X$ is a regular element in $\g$ and the subspace
$\C[X]\ell\sset V$ has an $X$-stable complement.
\end{definition}

Let $\rel\sset\X$ be the set of all relevant pairs $(X,\ell)\in\g\times\P$.
It is clear that $\rel$ is
a $G$-stable dense subset in $\X$ containg, in particular, all pairs $(X,\ell)$
such that $\C[X]\ell=V$.

Now, fix a  pair $(\fl,\Om)\in \scal$. We observe
that the pairs  $(X,\ell)$ in
the corresponding stratum $\X(\fl,\Om)$ 
are either {\em all} relevant or not,
i.e., either
$\X(\fl,\Om)\sset\rel$, in which case we call
the stratum $\X(\fl,\Om)$ {\em relevant}, or else
 $\X(\fl,\Om)\sset\X\sminus\rel$.

Explicitly, let $(X,\ell)\in \X(\fl,\Om)$ be such that
$X=z+x$ and $\ell=\C i$, where $z\in \fzr_\fl$ and
$(x,\ell)\in \Om$. Write
$\fl=\gl(V_1)\oplus\ldots\oplus \gl(V_l),\, X=\oplus_{k=1}^l\,X_k,$
and $i=i_1+\ldots+i_l,$
where $X_k\in \gl(V_k),\, i_k\in V_k$.
Then, the pair  $(X,\ell)$ is relevant
if and only if the following two conditions hold:

\pb{For each $k=1,\ldots,l$, in Jordan normal form, the element $X_k$
has  a single Jordan block.}

\pb{For any $k=1,\ldots,l,$ the vector $i_k$ is either
 a cyclic vector for
$X_k$ or else
$i_k=0$ (different alternatives for different values of
$k=1,\ldots,l$ are allowed).}

\subsection{The scheme $\La$.} \label{xreg}
Recall the  subscheme $\mnil\sset T^*(\g\times V)$
 defined by formula \eqref{nil}. We use  the identification in
\eqref{TP} and
 define a subscheme in $T^*\X$ as follows:
\begin{align}\label{Lambda}
\La:&=\left(\mnil\cap \bigl[T^*\g\times T^*(V\sminus\{0\})\bigr]\right)/
\C^\times\\
&=\{(X,Y,i,j)\in \g\times\g\times (V\sminus\{0\})\times V^*
\mid [X,Y]+ij=0\enspace\&\enspace Y\enspace \text{is nilpotent}\}/
\C^\times.\nonumber
\end{align}
 Note that the equation $[X,Y]+ij=0$  
implies that $\langle j,i\rangle=\Tr(ij)=
-\Tr[X,Y]=0$. 
Thus,  the  second line in \eqref{Lambda}
does define a subscheme of $T^*\X$.
It is clear that $\La$ is  a $G$-stable closed  subscheme of $T^*\X$.

The following result is a much more precise  version of
Theorem \ref{nilt}.

\begin{theorem}\label{mn} $\La$ is a Lagrangian subscheme
in $T^*\X$. More precisely,
with the notation of \S\ref{regstrata}, we have
$$\La\;=\;\bigcup_{\{(\fl,\Om)\in\scal
\;|\;\X(\fl,\Om)\sset \rel\}}\, 
\overline{T^*_{\X(\fl,\Om)}\X},
$$
where bar denotes the closure and the union on the right
is taken over all relevant strata.
\end{theorem}

\proof Fix a pair $(\fl,\Om)\in \scal,$
so, $\fl=\fz_\fl\oplus[\fl,\fl]$.
Let $p=(X,\ell)\in \X(\fl,\Om)$.
By $G$-equivariance,
it suffices to analyze the situation
in the case where $X=z+x$ and $\ell=\C i,$ for some $z\in \fzr_\fl$
and  $(x,i)\in \Om.$ 
We will  use the notation of \S\ref{regstrata}.
So, 
\beq{xz}
X=z+x,\en\text{where}\en 
x=\oplus_{k=1}^l\,x_k,\en x_k\in\gl(V_k)\en\text{is nilpotent, and}\en
z=\oplus_{k=1}^l\,z_k\cd\Id_{V_k},
\eeq for some
{\em pairwise distinct} complex numbers
$z_1,\ldots,z_l\in\C$.

\step{1.} We  need a more explicit description of
various  tangent and
cotangent
bundles.
The tangent space to $\X$ at $p=(X,\ell)$
has a direct sum decomposition
$T_p\X=\g\oplus T_\ell\P,$ where
$T_\ell\P$ denotes the tangent space to $\P$ at the point $\ell\in\P.$
Similarly, for cotangent spaces, we have  $T^*_p\X=\g^*\oplus T^*_\ell\P.$
We write $\BT_p\sset T_p\X $ for the
tangent, resp. $\BN_p\sset T^*_p\X$ for the
{\em conormal}, space to the stratum $\X(\fl,\Om)$
at $p$.

Let $S\sset \g$ be any subset.
It will be convenient to use shorthand notation and,
given a subset $U$ in either $\X=\g\times \P$
or in $T_p\X=\g\oplus T_\ell\P,$ resp.
in   $T^*_p\X=\g^*\oplus T^*_\ell\P$,
write $S+U:= \{(s+x',\ell')\mid s\in S,\,(x',\ell')\in U\}.$

With this notation,  we have $s+ L(x',\ell')=L(s+x',\ell'),$
 for any
$s\in  \fz_\fl,$ in particular, 
$\fzr_\fl+\Om=L(\fzr_\fl+x,\ell).$
Therefore, we get
$\X(\fl,\Om)=G(\fzr_\fl+\Om)=G\big(L(\fzr_\fl+x,\ell)\big).$
Hence,
$\X(\fl,\Om)=G(\fzr_\fl+x,\ell)$ and, for the tangent spaces, 
we deduce $\BT_p=\fz_\fl+ T_pG(p),$
where $T_pG(p)$ stands for the tangent space
at $p$ to the $G$-orbit through $p$.
Dually, we obtain
\beq{BN2}
\BN_p=(\fz_\fl\oplus 0)^\perp
\cap T^*_{G(p)}\X=(\fz_\fl^\perp \oplus T^*_\ell\P)\cap \mu\inv(0),
\eeq
where in the rightmost equality
we have used that the preimage of $0$ under the moment map is
the union  of the conormal bundles to the $G$-orbits in $\X$.
From now on, we will
 identify $\g^*$ with $\g$  via the trace form,
and thus view $\fz_\fl^\perp$ as a subspace in $\g$.

\step{2.} We begin  the proof of the Theorem by showing that
\begin{equation}\label{lan}
\La\sset
\bigsqcup_{(\fl,\Om)\in\scal}\, T^*_{\X(\fl,\Om)}\X.
\end{equation}

To this end, let $(X,Y,i,j)\in \La$. Thus,
$Y$ is nilpotent and $[X,Y]+ij=0.$ 
Hence, there exists a complete flag in $V$ that
is stable under the action of both $X$ and $Y$.
Therefore, it is also stable under the action
of $z$. Thus, all three matrices
$X,z,Y$ can be simultaneously made upper-triangular.
Furthermore, since $Y$ is nilpotent, in the
 upper-triangular form, $Y$ has
vanishing diagonal entries. We conclude that
$\Tr(z^m\cdot Y)=0$ for all $m=1,2,\ldots.$
Since $z=\oplus_{k=1}^l\,z_k\cd\Id_{V_k},$ cf. \eqref{xz},
we see by the  Vandermonde determinant,
that $\Tr(s\cdot Y)=0$ for any $s\in\fz_\fl.$
Thus, $Y\in \fz_\fl^\perp$ and \eqref{BN2} shows that
$(Y,j)\in\BN_p.$ This yields
\eqref{lan}.

\step{3.}
By Theorem \ref{nilt},  we
know that $\La$ is  a closed
Lagrangian scheme. This  
Lagrangian scheme is (set-theoretically) contained in
the RHS of \eqref{lan}, which is also  a Lagrangian scheme.
It follows that
 each irreducible component of
$\La$ must be at the same time an irreducible component of the
scheme   in
the RHS of \eqref{lan}, hence has the form 
$\overline{T^*_{\X(\fl,\Om)}\X},$ for
some pair $(\fl,\Om)$ such that
$T^*_{\X(\fl,\Om)}\X\sset\La$.
Thus, we have proved that  
$$\La\;=\;\bigcup_{\{(\fl,\Om)\in\scal\;|\;
 T^*_{\X(\fl,\Om)}\X\sset \La\}}\, 
\overline{T^*_{\X(\fl,\Om)}\X}.$$

We see that  completing the proof of the theorem
amounts to showing that the stratum
$\X(\fl,\Om)$ is relevant if and only if
$T^*_{\X(\fl,\Om)}\X\sset \La$, that is,
if and only if for some (hence any)
point $p=(X,\ell)\in\X(\fl,\Om),$ the following holds
\begin{equation}\label{condition}
(Y,j)\in\BN_p \quad\Longrightarrow
\quad Y\enspace\text{is  nilpotent}.
\end{equation}

\step{4.} We claim first that if $X$ is not regular in $\g$ then
\eqref{condition} does not hold.

We argue by contradiction.
If $X$ is not regular in $\g$ then
 there exists $k\in\{1,\ldots,r\}$ such that Jordan normal form
of the matrix $X_k$ consists of more than one
block. Hence, there exists a nonzero
semisimple element $Y_k\in \sv(V_k)$ that commutes
with $X_k$. 
We let $Y$ be the matrix whose
$k$-th component equals $Y_k$ and all other components
vanish.  Clearly,  for any
$s\in \fz_\fl,$ we have  $\Tr(s\cdot Y)=0,$
so $Y\in \fz_\fl^\perp$. 
Further, the quadruple $(X,Y,i,j)$,
where $i\in\ell$ and $j=0$ satisfies the
equation $[X,Y]+ij=0$. Thus, $(Y,j)\in\BN_p$.
Moreover, it is clear that $(X,Y,i,j)\not\in\La$ since
$Y$ is a nonzero semisimple, hence not nilpotent,
element.
Our claim  is proved.

Next, for $p=(X,\ell)$ as in \eqref{xz}, we claim
that
\beq{block_claim}
p=(X,\ell)\en\text{is relevant iff one has:}
\quad (Y,j)\in\BN_p\en\&\en
Y\in\fl\quad\Longrightarrow\quad Y\in \mn.
\eeq

We write $X=\oplus_{k=1}^l\, X_k$, where $X_k\in \gl(V_k)$. 
Clearly, 
it suffices to prove our claim for each $X_k$ separately.
Therefore, we may  assume that $\fl=\g$, so $X=z\cdot \Id_V+ x$,
for some $z\in\C$ and $x\in\mn$.
By the claim at the beginning of Step 4, we may further
restrict our attention to the
case where
$X$ is regular. Thus, we have reduced the proof of  \eqref{block_claim} to
the following

\step{5.}   Proof of \eqref{block_claim} assuming that $X$ consists of a single
Jordan $n\times n$-block. 

Recall that $\ell=\C\cdot i$,
and let $m\in \{0,1,\ldots,n\}$ be the unique integer such that
$i\in \Ker(x^m)\sminus\Ker(x^{m-1})$,
where we put $\Ker(x^0):=\{0\}$ and $\Ker(x^{-1}):=\emptyset$.
Let
$(Y,j)\in\BN_p$. We write $i=(i_1,\ldots,i_n)$ 
and $j=(j_1,\ldots,j_n)$, as usual. Lemma 
\ref{nak} shows that, conjugating the triple
$(Y,i,j)$ by an element of the group $G^X$
if necessary, we may assume that
\begin{equation}\label{ij}
i=(\underbrace{1,\ldots, 1}_{m},0,\ldots,0),\quad\text{and}
\quad
j=(\underbrace{0,\ldots,0}_{m}, j_{m+1},\ldots,j_n),
\end{equation}
for some (not necessarily nonzero)
$j_{m+1},\ldots,j_n\in\C$.

Now, writing out  the equation $[X,Y]+ij=0$
for  $i,j$ as above and $X$ an $n\times n$ Jordan block,
we find that
$Y$ must be upper triangular and, moreover, that
$$Y_{rs} = Y_{r-1, s-1} - i_{r-1}j_s,\quad\text{for all}\quad  1\leq r\leq s\leq n.$$
We see that one can choose, arbitrarily, the first row of $Y$
and then all other entries of $Y$ are uniquely determined
from the equations above. Solving these equations recursively
yields
\begin{equation}\label{yformula}
Y_{rs} = Y_{1, s-r+1} - i_{r-1}j_s - i_{r-2}j_{s-1} -\ldots-i_1j_{s-r+2},
\quad\text{for all}\quad  1<r\leq s\leq n.
\end{equation}
In particular, using \eqref{ij}, for $r=s$ we find
\begin{equation}\label{Y}
Y_{rr}=
\begin{cases}
Y_{11}&\text{if}\enspace 0< r\leq m \leq n \\
Y_{11}-j_{m+1}&\text{if}\enspace 0<m<r\leq n\\
Y_{11}&\text{if}\enspace 0=m<r\leq n.
\end{cases}
\end{equation}

Assume first that  $0<m<n$. We choose some complex number   $j_{m+1}\neq 0$,
and find $Y_{11}$ from the equation
$n\cdot Y_{11}-(n-m)\cdot j_{m+1}=0$. Then,  formula \eqref{Y} 
shows that if we let $j$ be any row vector with the chosen $(m+1)$-th
entry
$j_{m+1}$, then there exists a {\em trace zero}
upper triangular matrix $Y$ that satisfies  $[X,Y]+ij=0$
and such that $Y_{11}\neq 0.$ This matrix $Y$ is clearly
{\em not} nilpotent, hence, we have found
a pair $(Y,j)\in\BN_p$ with non-nilpotent
$Y$. Thus, condition
\eqref{condition} does not hold for $X$.
Hence, the conormal bundle $T^*_{\X(\fl,\Om)}\X$ is not contained in $\La$.

Assume now that either $m=0$ or $m=n$. Then, formula  \eqref{Y} 
combined with the requirement that $\Tr Y=0$ forces $Y_{rr}=0,$
for all $r=1,\ldots,n$. Hence, we have shown
that, for any $(Y,j)\in\BN_p$,
the matrix $Y$ is neccesarily nilpotent.
Thus, condition
\eqref{condition} holds  for~$X$.

Finally, we observe that in the case $m=0$ we have
$i=0$, while in the case $m=n$ the
vector $i=(1,\ldots, 1)$ is a cyclic
vector for $X$. 
Thus, $(X,\C i)$ is a relevant pair.

This completes the proof of \eqref{block_claim}.

\step{6} Claim \eqref{block_claim} implies, in particular,
that if \eqref{condition} holds, then
the pair $(X,\ell)$ is relevant. Now, let $p=(X,\ell)$
be any relevant pair. To complete the proof of the Theorem,
we must show that \eqref{condition} holds for $p$.

To this end, we use 
Definition \ref{distinguished} and write
$V=V^{1}\oplus V^{2},$ where $V^{1}=\C[X]\ell$ and
$V^{2}$ is an $X$-stable complement.
Accordingly,  we write an arbitrary element
 $Y\in \gl(V)$ in  block form,
$Y=\|Y^{rs}\|,$ where $Y^{rs}\in\Hom(V^r,V^s),\,  r,s\in \{1,2\}.$ 
The matrix $X$ is block-diagonal, so we have $X=X^{1}\oplus X^{2}$, where
$X^{r}=z^{r}+x^{r}\in\gl(V^{r}),\, r=1,2,$  are both regular.

For each  $r,s\in \{1,2\}$ we have the map
$$(\ad X)^{rs}: \Hom(V^r,V^s)\too\Hom(V^r,V^s),
\quad
Y^{rs}\mto X^{s}\ccirc Y^{rs}-Y^{rs}\ccirc X^{r}.
$$
This map has Jordan decomposition
$(\ad X)^{rs}=(\ad z)^{rs}+(\ad x)^{rs}.$
Since the eigenvalues
$z_1,\ldots,z_l,$ of $z,$ are
pairwise distinct, see \eqref{xz}, the maps $(\ad X)^{12}$ and
$(\ad X)^{21}$ are both invertible.

Now, fix a pair $(Y,j)\in \BN_p,$  as in \eqref{BN2}.
Write $j=j^{1}\oplus j^{2}$. Also, we have
 $\ell=\C i,$ and $i\in V^{1}$ is a cyclic vector for $X^{1}$.
Writing out the equation  $[X,Y]+ij=0$ 
 block-by-block,
we find 
\beq{block1}
(\ad X)^{rs}(Y^{rs})+ i^sj^r=0,\quad\forall r,s\in\{1,2\}.
\eeq
Since $i\in V^{1}$
and $X^{12}$ is an invertible map, equation  \eqref{block1} forces
$Y^{12}=0$.
Therefore, to show that $Y$ is nilpotent,
it suffices to show that both
$Y^{11}$ and $Y^{22}$ are nilpotent.
Hence, we may ignore the block $Y^{21}$.
Replacing $Y^{21}$ by zero does not affect equation
  \eqref{block1} for diagonal blocks,
so from now on we
assume that $Y= Y^{11}\oplus Y^{22}$.

From equation  \eqref{block1}
for the block corresponding to $\Hom(V^{2},V^{2}),$
using  that  $i\in V^{1}$, we deduce
$[X^{2},Y^{22}]=0$. For the block corresponding to $\Hom(V^{1},V^{1})$,
we use that $i$ is a cyclic vector of the operator
$X^{1}$ and
apply Lemma \ref{nak}. We deduce that
$j^{1}=0,$ hence equation  \eqref{block1}
yields $[X^{1},Y^{11}]=0.$ We see that $[X,Y]=0$.
This implies that $Y$ commutes with $z$, hence
we get $Y\in\fl.$

Thus, we are in the situation as at the end of Step 4.
Specifically, since $(X,\ell)$ is relevant,
from \eqref{block_claim} we obtain that
$(Y,j)\in\La$. This completes the proof of
the Theorem.
\endproof

\subsection{`Fourier dual' description of irreducible components.}
We recall the standard canonical isomorphism of symplectic manifolds:
\begin{equation}\label{TT}
T^*(\g\times\P)=\g\times\g^*\times T^*\P\stackrel{\,\g\leftrightarrow\g^*}\longeq
\g^*\times\g\times T^*\P=T^*(\g^*\times\P).
\end{equation}
Explicitly, using the identification $\g^*=\g$ and formula \eqref{TX} 
one rewrites isomorphism \eqref{TT} in  down-to-earth terms:
\begin{equation}\label{TXT}
T^*(\g\times \P)=\{(X,Y,i,j)\in \g\times\g\times (V\sminus \{0\})\times V^*
\mid \langle j,i\rangle=0\}\big/\C^\times
=T^*(\g^*\times \P),
\end{equation}
where the matrix $X$ is viewed as an element of $\g$
while  the matrix $Y$ is viewed as an element of~$\g^*$.

We have  vector bundle projections:
$$p:  T^*(\g\times\P)\to \g\times\P,\enspace
(X,Y,i,j)\mapsto (X,i);
\enspace
 p^\vee:  T^*(\g^*\times\P)\to \g^*\times\P,\enspace
(X,Y,i,j)\mapsto (Y,j).
$$

Let $\La^\vee\sset T^*(\g^*\times\P)$
be the image of $\La\sset T^*(\g\times\P)$ under
the canonical
isomorphism in \eqref{TT}.
It is clear that $\La^\vee$ is a closed
$G$-stable Lagrangian subscheme of $ T^*(\g^*\times\P)$.
Furthermore,  if
we view $\mn\times\P$ as a closed subset
in $\g^*\times\P$ (rather than in $\g\times\P$)
then, by definition of $\La$, we have
 $p^\vee(\La^\vee)\sset\mn\times\P$.

\begin{lemma}\label{La2}
The irreducible components of $\La^\vee$ are the sets
$\overline{T^*_S(\g^*\times\P)}$, where
$S$ runs through the (finite) set of all
$G$-diagonal orbits in $\mn\times\P$.
\end{lemma}

\proof It is clear that, for any smooth $G$-stable
locally closed subvariety $S$ contained
in $\mn\times\P$, we have
$T^*_S(\g^*\times\P)\sset \La^\vee$.
Thus, since $\La^\vee$ is closed,
we deduce that $\overline{T^*_S(\g^*\times\P)}\sset \La^\vee$,
for any $G$-diagonal orbit $S\sset\mn\times\P$.

To prove the converse, observe that $\La^\vee$
is  stable under dilations along the fibers of 
the projection $p^\vee:  T^*(\g^*\times\P)\to \g^*\times\P$,
i.e., it is a Lagrangian cone-subvariety. But any irreducible
closed Lagrangian cone-subvariety  $Z\sset T^*(\g^*\times\P)$
must be of the form $Z=\overline{T^*_S(\g^*\times\P)},$
where $S$ is an irreducible smooth $G$-stable
locally closed subvariety of $\g^*\times\P$, see e.g.
\cite{CG}, Lemma 1.3.27. Moreover, we 
have $S\sset p^\vee(Z)\sset \mn\times\P$. Therefore, we see from
Corollary \ref{mnfinite} that $S$ contains a unique
open dense $G$-diagonal orbit. Replacing $S$
by that orbit doesn't affect the closed set
$\overline{T^*_S(\g^*\times\P)},$
and the lemma follows.
\endproof

\section{A category of holonomic $\D$-modules}\label{cat}
\subsection{Twisted differential operators.}\label{tdo}
 For each $c\in\Z$, let $\oo_\P(c)$ be the corresponding  standard 
invertible sheaf on $\P$, and $\oo_\X(c)$ its pull-back via the second
projection $\X=\g\times\P\to\P.$ Write
$\D_\X(c)$ for the sheaf of algebraic {\em twisted differential operators}
 on $\X$
acting on the sections of $\oo_\X(c)$. Although the sheaf
 $\oo_\X(c)$ exists for integral values of $c$ only,
the corresponding  sheaf  $\D_\X(c)$ is well-defined for {\em any}
$c\in\C$, cf. \cite{BB2}. We write $\D(\X,c):=\Gamma(\X,\D_\X(c))$
for the algebra of global sections of the sheaf $\D_\X(c)$.

The action of any element of the Lie algebra $\g$ gives rise
to a vector field on $\g\times V$. Thus, we have
a morphism of the Lie algebra $\g$ into
the  Lie algebra of first order  differential operators
on $\g\times V$. The latter morphism extends uniquely
to a  filtration preserving
algebra homomorphism $\tau: \Ug\to\D(\g\times V).$
For each $c\in\C$, there is also a similar algebra homomorphism $\tau_c:
 \Ug\to\D(\X,c).$

Let $\eu$ denote the identity matrix viewed as
a base element in the  center of the Lie algebra $\g$.
It is clear that the $\ad\eu$-action  on $\g$ is trivial and
the action of $\eu$
 on $V$ generates the $\C^\times$-action on $V$ by  dilations.
It follows from definitions that
$\tau_c(\eu)=c,$ for any $c\in\C$. Furthermore, one can show that
the following canonical algebra morphisms are, in effect,
both isomorphisms
\beq{compare}
\D(\g\times V)^{\ad\eu}/\D(\g\times V)^{\ad\eu}\cd(\tau(\eu)-c)\iso
\Big(\D(\g\times V)/\D(\g\times V)\cd(\tau(\eu)-c)\Big)^{\ad\eu}\iso
\D(\X,c).
\eeq
Here, the first isomorphism is a special case
of \eqref{BA}, and the second isomorphism follows
from the
known description of the algebra $\G(\P,\D_\P(c))$.

\subsection{}
Let $\bz$ be the algebra of $\Ad G$-invariant {\em constant coefficient}
differential operators on $\g$. Thus, we have a natural
algebra isomorphism $\bz\cong(\Sym\g)^{\Ad G}.$
Write $\bz_+$ for the augmentation ideal in $\bz$ formed
by differential operators without constant term.
Now, any differential operator on $\g$ may be
identified with a differential operator on
$\g\times\P$ that acts trivially along the $\P$-factor.
This way, the obtain an algebra map
$\imath:\ \bz\into{\dx^G}.$

Let
 $\g=\C\cdot\eu\oplus\slv$ be an obvious
 Lie algebra direct sum decomposition. We
 put $\gc:=\tau_c(\slv)\sset \dx.$
Further,
let $\euu\in \D(\X,c)$ denote a first order differential
operator corresponding to the {\em Euler vector
field} along the factor $\g$ in the cartesian product
$\X=\g\times\P$. The differential operator $\euu$
is clearly $G$-invariant, i.e.,
$\euu$ commutes with $\gc$.
We write $\bu$ for the associative subalgebra in $\dx$ generated
by $\gc$ and $\euu$.

Motivated by \cite{Gi1} and \cite{Lu1}, \cite{Lu2},\cite{Lu3},
we are going to introduce a subcategory of the
category of finitely generated left $\dx$-modules.

\begin{definition}\label{gi} 
A finitely generated $\D(\X,c)$-module $M$ is called
{\em admissible} if the following holds:

\pb{The action  on $M$ of  the subalgebra $\imath(\bz_+)$ is {\em locally-nilpotent},
i.e., for any $u\in M$ there exists an integer $k=k(u)\gg 0$
such that $\imath(\bz_+)^ku=0$;}

\pb{The action  on $M$ of  the subalgebra
${\bu}$ is {\em locally-finite}.}
\vskip 1pt

\noindent
Let $\md$ be the full subcategory of $\Lmod{\dx}$
whose objects are admissible $\dx$-modules.
\end{definition}

\begin{remark}\label{G} In \cite{Gi3}, we study the group analog of the
notion of admissible
$\D$-module. The space $\X=\g\times\P$ is replaced there
by the space $\X^G:=G\times\P$. The group $G$ acts on
itself by conjugation and this makes $\X^G$
a $G$-variety with respect to diagonal action. This
way, one gets an algebra map $\tau^G_c: \U(\gc)\map
\D(\X^G,c)$. 

Further, let $\bz^G\cong(\Ug)^{\Ad G}$ be the algebra
of {\em bi-invariant} differential operators on
$G$ and $\bz^G_+:=\bz^G\cap(\g\cdot\Ug)$ its augmentation ideal.
One constructs similarly an algebra imbedding
$\imath^G: \bz^G\to \D(\X^G,c)$. Thus, the  two
conditions of Definition \ref{gi}  may be
replaced by their $G$-counterparts,
with $\bz_+$ being replaced by $\bz^G_+$,
resp.  $\bu$ being replaced by $\U(\gc)$. This
leads to the notion of  admissible
$\D(\X^G,c)$-module (the extra condition of local
finiteness of $\euu$-action has no group analogue
and may be dropped from definition).
\hfill$\lozenge$\end{remark}

\begin{remark}\label{SLint}
Recall that a (possibly infinite dimensional) $\gc$-module
 $M$ is said
to be {\em locally-finite} if $\dim(\U(\gc)\cdot m)<\infty$, for any $m\in M$.

It follows from  complete reducibility
of finite dimensional $\slv$-modules
that any
locally finite $\gc$-module splits into
(possibly infinite) direct sum of finite
 dimensional simple $\slv$-modules.
Furthermore,
 the  group $SL(V)$ being {\em simply-connected},
any finite dimensional $\slv$-module can be exponentiated
to  a  rational representation of the algebraic group $SL(V)$.
Thus, a locally finite $\gc$-module
is the same thing as a (possibly infinite) direct sum of finite
 dimensional simple rational representations
of the  algebraic group $SL(V)$ (on which the central
element $\eu$ acts as~$c\cdot\Id$).
\end{remark}

\subsection{}
The projective space $\P$ is a 
partial flag manifold for $G$,
hence, Beilinson-Bernstein theorem \cite{BB1} holds for $\P$.
The space $\g$ being affine, from \cite{BB1} we deduce
\begin{proposition}\label{BB} For any $c\geq 0$,
the functor of global sections provides an equivalence
between the abelian categories of $\D_\X(c)$-coherent sheaves
and of finitely-generated $\D(\X,c)$-modules, respectively.\qed
\end {proposition}

Given a  $\D_\X(c)$-coherent sheaf $\lll$, we write
$\supp\lll$ for its {\em support} as an $\oo_\X$-module,
and
$\CC\lll$ for its {\em characteristic cycle}, a  $\C^\times$-stable 
algebraic cycle in $T^*\X$. Abusing the notation, we will also
write $\CC\lll$ for the set-theoretic  support of the algebraic cycle 
$\CC\lll$.
Thus, we have $p(\CC\lll)=\supp\lll$,
where $p: T^*\X\to\X$ is the projection.

Using Proposition \ref{BB}, we will often identify an
object of the category $\md$ with a $\D_\X(c)$-coherent sheaf.
Thus, we write $\CC{L}$ for the
characteristic cycle, resp., $\supp{L}$ for the support,
 of the  $\D_\X(c)$-coherent sheaf
corresponding to an object $L\in \md$.

\begin{proposition}\label{SS} 
Let $M$ be a 
finitely generated $\dx$-module which is locally finite
as a ${\bu}$-module. Then,  we have
$M\in\md$ if and only if
$\CC{M}\sset\La$.

Furthermore, any object of $\md$ is a regular holonomic (twisted)
$\D$-module which is smooth along the strata $\X(\fl,\Om),
\,(\fl,\Om)\in\scal$.
\end{proposition}

\proof 
We will exploit the theory of  Fourier transform of $\D$-modules, 
see \cite{Br}.

Recall that the Fourier transform gives an equivalence
\begin{equation}\label{ff}
\ff_\D:\ \Lmod{\D(\g\times\P,c)}\iso\Lmod{\D(\g^*\times\P,c)},
\quad M\mto\ff_\D M,
\end{equation}
from the category of finitely generated (twisted)  $\D$-modules on  the
total space of the
(trivial) vector bundle $\X=\g\times\P\to\P$
to a similar category of  finitely generated twisted  $\D$-modules on 
the total space of the dual vector bundle $\g^*\times\P\to\P$. 

We will generalize earlier notation slightly,
and write $\euu$ for the Euler vector field
along the fibers of any vector bundle,
e.g., along the fibers of the vector bundle $\g\times\P\to\P$
or of  $\g^*\times\P\to\P$. Recall that a (twisted)
$\D$-module on a vector bundle is said to be
{\em monodromic} if it is locally finite as a $\euu$-module.

It is known that the Fourier transform functor 
acts especially nicely on monodromic modules.
In particular, it takes monodromic  (twisted) $\D$-modules
into  monodromic  (twisted) $\D$-modules.
Furthermore,  it is known (see \cite{Br}) that
given  a monodromic  (twisted) $\D$-module $M$, one has:
\smallskip

\pb{$\CC(M)=\CC(\ff_\D M)$ under the canonical isomorphism in \eqref{TT}.}

\pb{$M$ has  regular singularities
$\quad\Longleftrightarrow\quad
\ff_\D M$ has  regular singularities.}
\smallskip

From now on, we fix  an $\bu$-locally finite, finitely
generated $\dx$-module $M$
 such that $\CC(M)\sset\La$, and
let  $(\Sym\g)^{\Ad G}_+$ be the set of $G$-invariant polynomials
on $\g^*\cong\g$ without constant term.

It is immediate
from the  equation $\CC(M)=\CC(\ff_\D M)$ and definition of $\La$
that we have $\supp(\ff_\D M)\sset \mn\times\P$.
Any element of $(\Sym\g)^{\Ad G}_+$
vanishes on $\mn$.
Hence,  the action of $(\Sym\g)^{\Ad G}_+$ on
$\ff_\D M$ is locally nilpotent.
It follows that the action of the subalgebra
$\imath(\bz_+)$ on $M$ is also  locally nilpotent.
This proves the implication
$\CC{M}\sset\La\enspace\Longrightarrow\enspace
M\in\md$. 

Next, we have

\begin{claim}\label{claim}
Any $\Ug$-locally finite, finitely generated $\dx$-module 
$M$ such that $\supp M\sset \mn\times\P$
is automatically a holonomic
module with regular singularities.
\end{claim}

\proof[Proof of Claim.]
As we have explained in Remark \ref{SLint}, any
$\U(\gc)$-locally finite $\dx$-module is automatically
$SL(V)$-{\em equivariant}.
It is clear that
any  $G$-diagonal orbit
 $S\sset \mn\times\P$ is also
a single $SL(V)$-orbit. 
It follows that the algebraic variety   $\mn\times\P$ is partitioned
into finitely many $SL(V)$-orbits, by Corollary \ref{mnfinite}.
But it is a well-known result that if $H$ is a
linear algebraic group and ${\mathscr{X}}$ is an
arbitrary $H$-variety
with finitely many $H$-orbits then,
any $H$-equivariant $\D$-module on
${\mathscr{X}}$ is a holonomic
module with regular singularities.
This proves that, for any $\bu$-locally finite, finitely generated
$\dx$-module $M$ such that $\CC(M)\sset \La$,
the $\D$-module $\ff_\D M$ has regular singularities.
Therefore, by the properties of Fourier transform 
mentioned above,  $M$  also has regular singularities.
This proves Claim \ref{claim}.
\endproof

To complete the proof of the Proposition, we must 
show that
$M\in\md\enspace\Longrightarrow\enspace\CC{M}\sset\La$.
This
is entirely analogous to the proof
of implication (i) $\Longrightarrow$ (ii) in
 \cite{Gi1}, Theorem 1.4.2. We leave details to the reader.
\endproof

\begin{corollary} $\md$ is an abelian, artinian category
with finitely many simple objects; in particular,
 every object of $\md$ has finite length.\qed
\end{corollary}

\begin{remark} There is a natural group analog of
the Lagrangian subscheme $\La\sset T^*(\g\times\P)$. It is a 
Lagrangian subscheme $\La^G\sset T^*(G\times\P)$ defined by
$$\La^G=\{(g,Y,i,j)\in G\times\g\times (V\sminus\{0\})\times V^*
\mid \Ad g(Y) -Y+ij=0\enspace\&\enspace Y\enspace \text{is nilpotent}\}/
\C^\times.
$$
One can then formulate and prove a group version of
Proposition \ref{SS} for $\D(\X^G,c)$-modules, cf.
Remark \ref{G}.
\end{remark}

\section{Cherednik algebra and Hamiltonian reduction.}
\subsection{Reminder.}\label{reminder} Recall
 the permutation representation  $\h=\C^n$ of $W=S_n$, the Symmetric group.
We write $y_1$, $\ldots$, $y_n$ for the standard basis of $\h,$
and $x_1$, $\ldots$, $x_n\in\h^*$ for the corresponding
dual coordinate functions.
Let $s_{ij}\in S_n$ denote the 
transposition $i\leftrightarrow j$.

Let $c\in\C$. The 
\emph{rational Cherednik algebra} of type $\mathbf{A}_{n-1}$,
as defined in \cite{EG}, is an associative $\C$-algebra $\hh_c$ with generators
$x_1$, $\ldots$, $x_n$, $y_1$, $\ldots$, $y_n$ and
the group $S_n$, and the following defining relations:
\begin{equation}\label{defrel}
\begin{array}{lll}
\displaystyle
&{}_{_{\vphantom{x}}}s_{ij}\cdot x_i=x_j\cdot s_{ij}\quad,\quad
 s_{ij}\cdot y_i=y_j\cdot s_{ij}\,,&
\forall i,j\in\{1,2,\ldots,n\}\;,\;i\neq j\break\medskip\\
&{}_{_{\vphantom{x}}}{}^{^{\vphantom{x}}}
[y_i,x_j]= c\cdot s_{ij}\quad,\quad[x_i,x_j]=0=[y_i,y_j]
\,,&
\forall i,j\in\{1,2,\ldots,n\}\;,\;i\neq j\break\medskip\\
&{}^{^{\vphantom{x}}} [y_k,x_k]= 1-c\cdot\sum_{i\neq 
k}\;s_{ik}\,.
\end{array}
\end{equation}

In \S\ref{hhham}, we have mentioned
the spherical  subalgebra $\ehe\subset \hh_c$.
The assignment $a\mapsto a\cdot\e=$\break
$\e\cdot a$
gives algebra imbeddings
$(\Sym\h)^W=\C[y_1,\ldots,y_n]^W\into \ehe$
and $\C[\h]^W=\C[x_1,\ldots,x_n]^W$
$\into\ehe$.
We will identify the algebras
$(\Sym\h)^W$ and $\C[\h]^W$ with their images in $\ehe$.
It is known that these two algebras generate
$\ehe$ as an algebra.

It is immediate to see from the defining relations
\eqref{defrel}
 that the following assignment 
\beq{fourier}
w\mto w\en (\forall w\in S_n),\quad x_i\mto y_i, \quad y_i\mto -x_i,
\quad \forall  i=1,\ldots,n,
\eeq
extends to an algebra automorphism
$\hh_c\to\hh_c$, called {\em Fourier  automorphism}, cf. \cite[\S7]{EG}.
By restriction, we also get the 
automorphism $\ehe\to\ehe$ and, by transport of structure, an auto-equivalence
$\ff_\hh: \Lmod{\ehe}\to\Lmod{\ehe},$
called Fourier transform functor.

The Cherednik algebra has an 
 increasing filtration such that
all elements of $S_n$ as well as the generators
 $x_1$, $\ldots$, $x_n\in \hh_c$ have
filtration degree zero, and the generators
 $y_1$, $\ldots$, $y_n\in\hh_c$ have
filtration degree $1$. We equip $\ehe$ with the
induced filtration. Then, the Poincar\'e-Birkhoff-Witt
theorem for Cherednik algebras, see \cite{EG},
yields a graded algebra isomorphism
\begin{equation}\label{PBW}
\gr(\ehe)\cong\C[\h\times\h^*]^W.
\end{equation}

The rest of this section is devoted to the
proof of our main result,  Theorem \ref{iso}.

\subsection{}\label{preliminary1}  We have the
moment map $\mu_{_\X}: T^*\X\to\g^*$
and its scheme-theoretic zero fiber:
$$\mu_{_\X}\inv(0)=\{(X,Y,i,j)\in \g\times\g\times (V\sminus \{0\})\times V^*
\mid [X,Y]+ij=0\}\big/\C^\times.
$$
(as has been explained in \S\ref{xreg}, this is indeed a subscheme
in $T^*\X$ since the equation
$[X,Y]+ij=0$ implies $\langle j,i\rangle=0$).
It is clear that $\mu_{_\X}\inv(0)$ is a closed $G$-stable subscheme 
of $T^*\X$, furthermore,  Theorem \ref{t1}(iii)
implies that  $\mu_{_\X}\inv(0)$ is a reduced complete intersection.
By definition, restriction of functions from
$T^*\X$ to $\mu_{_\X}\inv(0)$ gives an algebra
isomorphism:
\beq{mu0}
\C[T^*\X]/\C[T^*\X]\cdot\mu_{_\X}^*(\gc)=
\C[\mu_{_\X}\inv(0)],
\eeq
where the Lie algebra $\gc$ is
identified with the vector space of linear functions on $\g_c^*$.

\begin{lemma}\label{iso_lemma}
There is a graded algebra isomorphism
$$
\C[\mu_{_\X}\inv(0)]^G
=\bigl(\C[T^*\X]/\C[T^*\X]\cdot\mu_{_\X}^*(\gc)\bigr)^G\cong \C[\h\times\h]^W.$$
\end{lemma}
\proof We write $\OO_1\sset\g$ for the conjugacy class of
rank one nilpotent matrices, and let $\OOO=\OO_1\sqcup\{0\}$
be the 
 closure of $\OO_1$ in $\g$. The moment map
$\pi: T^*\P\to\g^*=\g,\, (i,j)\mapsto ij$,
see formula \eqref{TP},  gives
a birational isomorphism $T^*\P\map\OOO$.

Let $\I_1\sset\C[\g]=\C[Z]$ be the ideal generated
 by all the $2\times2$ minors of the matrix $Z$ and also by the 
linear function $Z\mapsto \Tr Z.$
This is  known to be a prime ideal whose zero scheme
equals  $\OOO$.
Furthermore, the  pull-back morphism
$\pi^*: \C[\g]/\I_1=\C[\OOO]\to \C[T^*\P]$
is   known to be a graded algebra isomorphism.

Next, we write $T^*\X=T^*\g\times T^*\P$.
It is clear that the
moment map $\mu_{_\X}: T^*\X\to\g^*$ may be factored as
the following composite map:
\beq{factored}
T^*\X=T^*\g\times T^*\P
\stackrel{\pi_{_\X}}\tooo
\g\times\g^*\times\OOO=\g\times\g\times\OOO
\stackrel{\theta}\too\g^*=\g,
\eeq
where $\pi_{_\X}:=\Id_{T^*\g}\times\pi,$ and
where the map $\theta$ is given by
$\g\times\g\times\OOO\ni (X,Y,Z)\mto [X,Y]+Z$.
Observe that, by the last sentence of the
preceeding paragraph, the map  $\pi_{_\X}$
induces  a graded algebra isomorphism
\beq{pix}\pi_{_\X}^*:\
 \C[\g\times\g]\,\otimes\, (\C[\g]/\I_1)=
\C[\g\times\g^*\times\OOO]\iso \C[T^*\X].
\eeq

Write $\C[X,Y,Z]:=\C[\g\times\g\times\g]$,
and let $\C[X,Y,Z]([X,Y]+Z)$
 denote the ideal in $\C[X,Y,Z]$
generated by all  matrix entries of
the matrix $[X,Y]+Z$. We also consider a larger ideal
$\II:=\C[X,Y]\otimes \I_1+$
$\C[X,Y,Z]([X,Y]+Z)\sset\C[X,Y,Z].$
Thus, from \eqref{factored} and \eqref{pix}
we find  
$$
\C[T^*\X]/\C[T^*\X]\cdot\mu_{_\X}^*(\g)\cong
\C[\g\times\g\times\OOO]/\C[\g\times\g\times\OOO]
\cdot\theta^*(\g)=\C[X,Y,Z]/\II.
$$

We define an algebra homomorphism
$r: \C[X,Y,Z]\to \C[X,Y]$ to be 
the map sending a polynomial
$P\in 
\C[X,Y,Z]$ to the function
$(X,Y)\mapsto P(X,Y,-[X,Y])$.
The homomorphism $r$
  clearly induces an isomorphism
$\C[X,Y,Z]/\C[X,Y,Z]([X,Y]+Z)\iso\C[X,Y]$.
Observe that the linear function $P: (X,Y,Z)\mapsto \Tr Z=
\Tr([X,Y]+Z)$ 
belongs to the ideal $\C[X,Y,Z]([X,Y]+Z)$.
We deduce that the map $r$ sends
 the subspace
$\C[X,Y]\otimes \I_1\sset\C[X,Y,Z]$
 to $I_1\sset\C[X,Y]$,
the ideal $I_1$ considered in  Section \ref{around}.
Thus, we obtain  algebra isomorphisms
$$
\C[T^*\X]/\C[T^*\X]\cdot\mu_{_\X}^*(\gc)\cong\C[X,Y,Z]/\II\cong
\C[X,Y]/I_1.
$$

The statement of the Lemma is now  immediate
from isomorphisms  \eqref{composite}.
\endproof

\subsection{The Dunkl  homomorphism.}\label{dunkl}
Let $\hreg\sset \h=\C^n$ be the
Zariski open dense subset  formed by $n$-tuples with
pairwise distinct coordintates.
The Dunkl homomorphism is an algebra imbedding
$\Theta: \ehe\into \D(\hreg)^W$
defined as follows, see e.g. \cite{EG}.

Write $\C[\h^*]\# W$ for the smash product algebra,
and let $\triv: \C[\h^*]\# W\to\C$ be the homomorphism
that sends every element $w\in W$ to $1,$ and 
acts on the polynomial algebra $\C[\h^*]$ by $f\mapsto f(0).$
We 
view  $\C[\h^*]\# W$ as a subalgebra of $\hh_c$
and let $\hh_c\bigotimes_{\C[\h^*]\# W}\triv$
be the induced left $\hh_c$-module.

The natural imbedding
$\C[\h]\into\hh_c$ yields,
 by the Poincar\`e-Birkhoff-Witt theorem
for the Cherednik algebra $\hh_c,$ cf. \cite{EG},
a vector space isomorphism
$\C[\h]\iso 
\hh_c\otimes_{\C[\h^*]\# W}\triv$. The left $\hh_c$-action on 
$\hh_c\otimes_{\C[\h^*]\# W}\triv$ gets transported,
via the isomorphism, to an $\hh_c$-action on $\C[\h].$ 
The resulting action of the spherical
subalgebra $\ehe$ preserves the subspace
$\C[\h]^{W}=\e\cdot\C[\h],$ of symmetric polynomials. 
Moreover,
a direct
calculation
shows that,
for any $u\in \ehe,$ the corresponding action-map 
$u: \C[\h]^{W}\to \C[\h]^{W}$
is given
by a $W$-invariant differential operator, $\Th(u),$
 with rational coefficients. More precisely,
all coefficients of the  differential operator $\Th(u)$ turn
out to be  regular functions on $ \hreg$.

The assignment $u\mapsto\Th(u)$ gives the
desired  Dunkl homomorphism
$\Theta: \ehe\to \D(\hreg)^W$.
This is a filtration preserving injective algebra homomorphism,
and we let $\gr\Th$ denote the corresponding associated graded map.
The composite map
$$\C[\h\times\h^*]^W=\gr(\ehe)\stackrel{\gr\Th}\tooo
\gr(\D(\hreg)^W)=\C[\hreg\times\h^*]^W$$
is known to be equal to the natural
restriction map $(\jmath\times\Id)^*: \C[\h\times\h^*]^W\to
\C[\hreg\times\h^*]^W$,
induced by the open imbedding
$\jmath: \hreg\into\h=\C^n$.

\subsection{The radial part map.}\label{rad_sec} Recall that
$\gc$ denotes the image in $\dx$ of the Lie subalgebra
$\slv\sset \g$,
and  let  $\D(\X,c)\cdot{\gc}\sset \D(\X,c)$
be  the left ideal generated by $\gc$.

Using \eqref{compare}, for any $c\in\C,$ one obtains 
\begin{equation}\label{inv}
\bigl(\dx/\D(\X,c)\cdot{\gc}\bigr)^{\ad\gc}\cong
\Bigl(\D(\g\times V)\big/\D(\g\times V)\cd(\tau(\eu)-c)+ 
\D(\g\times V)\cd\tau(\slv)\Bigr)^{\ad\g}.
\end{equation}
Here, $\ad\g$-invariants are taken with respect to
the
`adjoint' action defined by the formula
$\ad g: u\mto \tau(g)\cdot u-u\cdot\tau(g),\,\forall g\in\g$
(or a similar formula with $\tau_c$ instead of $\tau$).
The object on each side of \eqref{inv}
is the result of a certain Hamiltonian reduction,
 see \S\ref{gen} below and \cite[\S3.4]{BFG}. 
Thus, each side in  \eqref{inv}
acquires a natural  algebra structure and
the isomorphism in  \eqref{inv} is an algebra isomorphism.
In this paper, we will only work with the algebra on the left hand side
of  \eqref{inv}. Thus, we will neither use nor prove the isomorphism in
  \eqref{inv}.

We now recall the construction, due to  \cite[Proposition 5.3.6]{BFG}, 
of the following 
 filtered algebra homomorphism, called
 `radial part' map
\begin{align}\label{rad_map}
\Psi_c:\ & (\dx/\dx\cd\gc)^{\ad\gc}
\map  \D(\hreg)^W,\quad\text{such that}\\
&\Theta(\C[\h]^W)=\Psi_c(\C[\g]^G),
\quad\text{and}\quad
\Theta\bigl((\Sym\h)^W\bigr)=\Psi_c(\bz).\nonumber
\end{align}

Let $\eer\sset\V=\g\times V$ be
the subset formed by
the pairs $(x,v)$ such that $v$ is a 
{\em cyclic vector} for the operator $x: V\to V$.
It is clear that $\eer$
a $G$-stable, Zariski open dense
subset of $\V$.
We compose  the first projection
$\V=\g\times V\to\g$ with the adjoint
quotient map $\g\to\g/\Ad G=\h/W$,
and restrict the resulting morphism
to the subset $\eer\sset\V$. This way we get
a morphism $p: \eer\to \h/W$.
It turns out that the group $G$
acts freely along the fibers of $p$, and this
makes  the map
$p: \eer\to \h/W$ a
principal $G$-bundle over $\h/W$,
cf. \cite[Lemma 5.3.3]{BFG}.

Let $\xreg$ be the image of $\eer$
under the projection $\g\times(V\sminus\{0\})\onto \g\times \P$.
The map $p$ descends to $\xreg$ and makes it  a 
principal $PGL(V)$-bundle on $\h/W$.
Now, the standard description of differential operators
on the base of a principal bundle in terms of those on the
total space of the bundle yields
an algebra isomorphism
\beq{rad2}
\big(\D(\xreg,c)/\D(\xreg,c)\cd \gc\big)^G\iso\D(\h/W),\quad\forall
c\in\C.
\eeq
 
We have  (strict) inclusions
$\D(\h)^W\sset\D(\h/W)\sset\D(\hreg)^W$.
The map $\Psi_c$ in \eqref{rad_map}  is defined as the following composite
homomorphism
$$\xymatrix{
\big(\dx/\dx\cd\gc\big)^G
\ar[rr]^<>(0.5){{_\text{restriction}}}&&
\big(\D(\xreg,c)/\D(\xreg,c)\cd \gc\big)^G\;
\ar[r]^<>(0.5){\eqref{rad2}}&
\;\D(\h/W)
}\into\D(\hreg)^W.
$$

More explicitly,
fix a nonzero volume element $\bv^*\in\wedge^n V^*$ and
let
$$
(x,v)\mto \ffs(x,v):=\langle\bv^*, v\wedge x(v)\wedge\ldots
\wedge x^{n-1}(v)\rangle,
$$
a polynomial 
function on $\V$.
Further, for any integer $c\in\Z$, we put
$$\oo(\eer,c):=\{f\in \C[\eer]\mid g^*(f)=(\det g)^c\cdot f,
\enspace\forall g\in G\}.
$$
Observe that  $\ffs\in \oo(\eer,1)$ and 
 $\eer=\V\sminus\ffs\inv(0)$,
in particular, $\eer$ is an affine
variety.

It is clear that pull-back via the bundle projection $p: \eer\to\h/W$
makes the vector space $\oo(\eer,c)$ a $\C[\h/W]$-module. Furthermore, one shows that
this is in effect a rank one free  $\C[\h/W]$-module
with generator $\ffs^c$, so one has
a bijection $\C[\h/W]\iso \oo(\eer,c),\, f\mto \ffs^c\cdot p^*(f)$.
The isomorphism in \eqref{rad2} is obtained by transporting
the action of differential operators on $\eer$
via the 
bijection.

\subsection{Proof of Theorem \ref{iso}.}\label{pf}
Our argument follows the strategy of 
\cite[\S5]{BFG} (that, in its turn, is based on an
argument from \cite{EG}), and we will freely use the notation
from~\cite{BFG}.
 
\step{1.} We claim  that the image of $\Theta$ is contained in
the image of $\Psi_c$, i.e., we have
\beq{sumup}
\Theta(\ehe)\subseteq\Psi_c\bigl((\dx/\dx\cdot\gc)^{\ad\gc}\bigr).
\eeq
To see this,  recall that the two subalgebras $\C[\h]^W,
\,(\Sym\h)^W\sset\ehe$
generate $\ehe$ as an algebra. It follows that
$\Theta(\C[\h]^W)$ and 
$\Theta\bigl((\Sym\h)^W\bigr)$
generate $\Theta(\ehe)$ as an algebra.
Thus, equations \eqref{rad_map} yield the inclusion in
\eqref{sumup}.

From \eqref{sumup} we deduce an imbedding
of the corresponding  associated graded algebras:
\beq{include}
\gr\Theta(\ehe)\,\into\,\gr\Psi_c\bigl((\dx/\dx\cdot\gc)^{\ad\gc}\bigr),
\eeq

\step{2.} For any smooth manifold $\mathcal{Y},$
the  sheaf $\D_{\mathcal{Y}}$, of differential (or twisted  differential)
operators on $\mathcal{Y},$  comes equipped with the
 standard increasing
filtration by the order of differential operator.
Let $\D({\mathcal{Y}})=\G({\mathcal{Y}},\D_{\mathcal{Y}})$
and write
$\gr\D({\mathcal{Y}})$ for the associated graded
algebra.
The
principal symbol  map provides
a canonical graded algebra imbedding
$\gr\D({\mathcal{Y}})\into \C[T^*{\mathcal{Y}}]$.

In the special case $\mathcal{Y}=\g$, the imbedding
 $\gr\D(\g)\into\C[\g\times\g^*]$ is clearly an isomorphism.
Also, for $\mathcal{Y}:=\P=\P(V)$,
the principal symbol  map $\gr\D(\P,c)\into\C[T^*\P]$ is well known to be
an isomorphism. We deduce that for $\X=\g\times\P$
the principal symbol  map yields an isomorphism
$\gr\D(\X,c)\iso\C[T^*\X].$

The imbedding $\gc\into \dx$ extends, by multiplicativity,
to a  filtered algebra map $\U(\gc)\to \dx$.
This  algebra map induces
 an associated graded  homomorphism 
$$\Sym\gc=\gr\U(\gc)\map\gr\D(\X,c)=
\C[T^*\X].$$
The latter homomorphism 
is well known to be
the pull-back via the moment map $\mu_{_\X}:T^*\X\to\g_c^*.$
Hence, in $\C[T^*\X]=\gr\D(\X,c)$, we have 
$\C[T^*\X]\cdot\mu_{_\X}^*(\gc)=
(\gr\D(\X,c))\cdot\gc\subseteq\gr\bigl(\D(\X,c)\cdot\gc\bigr).$
Thus, using \eqref{mu0}, we obtain the following chain of
graded algebra morphisms
$$
\C[\mu_{_\X}\inv(0)]=\C[T^*\X]/\C[T^*\X]\cd\mu_{_\X}^*(\gc)=
\gr\D(\X,c)/(\gr\D(\X,c))\cd\gc
\stackrel{_\text{proj}}{-\!-\!\!\!\twoheadrightarrow}
\gr\bigl(\D(\X,c)/\D(\X,c)\cd\gc\bigr).
$$

Restricting these morphisms  to
$\ad\gc$-invariants and using \eqref{BA},
we obtain the following graded algebra morphisms
\beq{wellG}
\C[\mu_{_\X}\inv(0)]^{\ad\gc}=
\Big(\gr\D(\X,c)/(\gr\D(\X,c))\cd\gc\Big)^{\ad \gc}
\stackrel{_\text{proj}}{-\!-\!\!\!\twoheadrightarrow}
\gr\Bigl(\bigl(\D(\X,c)/\D(\X,c)\cd\gc\bigr)^{\ad \gc}\Bigr).
\eeq

\step{3.} Let $\Id_\h\times\kappa:
\C[\h\times\h^*]^W\iso\C[\h\times\h]^W$
be  the algebra isomorphism
arising from the bijection
$\kappa: \h^*\iso \h$ induced by
 the trace pairing on $\g$. 
Thus,
we obtain the following diagram:
$$
\xymatrix{
\C[\h\times\h^*]^W\ar@{=}[d]^<>(0.5){\eqref{PBW}}
\ar@{=}[r]^<>(0.5){\Id_\h\times\kappa}
&\C[\h\times\h]^W\ar@{=}[r]^<>(0.5){\text{Lemma \ref{iso_lemma}}}&
\bigl(\C[T^*\X]/\C[T^*\X]\cd\mu_{_\X}^*(\gc)\bigr)^G
\ar@{=}[d]_<>(0.5){\eqref{wellG}}\\
\gr(\ehe)\ar@{=}[d]^<>(0.5){\Theta}&&
\left(\gr\dx\big/\bigl(\gr\dx\bigr)\cd\gc\right)^{\ad\gc}
\ar@{->>}[d]_<>(0.5){\text{proj}}\\
\gr\Theta(\ehe)
\ar@{^{(}->}[r]^<>(0.5){\eqref{include}}&
\gr\Psi_c\bigl((\dx/\dx\cdot\gc)^{\ad\gc}\bigr)&
\gr\left(\bigl(\D(\X,c)/\D(\X,c)\cd\gc\bigr)^{\ad \gc}\right)
\ar@{->>}[l]_<>(0.4){\gr\Psi_c}
}
$$

Now, it is straightforward to verify that this
 diagram of graded algebra maps {\em commutes}.
This forces both  surjections in the diagram,
as well as the injective  map \eqref{include} in the 
bottom row of the diagram, all to be {\em bijective}.
Therefore, we deduce the following isomorphisms:

\begin{equation}\label{equality}
\gr\bigl(\dx/\dx\cdot\gc\bigr)^{\ad\gc}\cong
\gr\Psi_c\bigl((\dx/\dx\cdot\gc)^{\ad\gc}\bigr)
\stackrel{\psi}=\gr\Theta(\ehe)\cong\gr(\ehe),
\end{equation}
where  the  equality indicated
as $\psi$ holds inside the bigger algebra
$\gr\D(\hreg)^W$.

As has been  proved earlier, the algebra
 $\Theta(\ehe)$ is contained in
$\Psi_c\bigl({(\dx/\dx\cdot\gc)^{\ad\gc}}\bigr).$
The equality $\psi$ may be obtained as the
associated graded map corresponding to the
 imbedding $\Theta(\ehe)\into$
$\Psi_c\bigl({(\dx/\dx\cdot\gc)^{\ad\gc}}\bigr).$
It follows that the imbedding itself is, in effect,
an equality (that holds in the larger algebra $\D(\hreg)^W$). 
Hence, we may invert the (injective) map
$\Theta$ and define a graded algebra morphism $\Phi_c$
as the following composite map, very similar to the
one used in \cite{EG}: 
$$\Phi_c: \
(\dx/\dx\cd\gc)^{\ad\gc}\stackrel{\Psi_c}\too
\Psi_c\bigl((\dx/\dx\cd\gc)^{\ad\gc}\bigr)=
\Theta(\ehe)\stackrel{\Theta\inv}\too
\ehe.
$$
It is immediate from \eqref{equality} that the corresponding
associated graded map gives a bijection
$\gr\Phi_c:\gr\bigl(\dx/\dx\cd\gc\bigr)^{\ad\gc}\iso
\gr(\ehe).$ Thus, the map $\Phi_c$ is itself
a bijection. The theorem is proved.\qed

\begin{corollary}\label{corG} The projection $\operatorname{proj}$
in \eqref{wellG} is a bijection; in particular, one has
a graded algebra isomorphism
 $\gr\left(\bigl(\dx/\dx\cdot\gc\bigr)^{\ad\gc}\right)\cong 
\C[\mu_{_\X}\inv(0)]^G$.\qed
\end{corollary}

\subsection{An application to \cite{EG}.}
We fix $c\in\C$ and use other notation of previous sections.
The natural $GL(V)$-action on $\P=\P(V)$ gives an
algebra homomorphism $\Ug\to\D(\P,c)$. The kernel
of this homomorphism is known to be a {\em primitive}
ideal $\ind\sset \Ug$,
moreover, it is exactly the
primitive ideal considered in \cite{EG}.

The group $G=GL(V)$ acts on $\g$ via
the adjoint action.
Differentiating this action gives rise to
an associative algebra homomorphism
$\ad: \Ug \to \D(\g)$. Let
$\D(\g)\cdot\ad\ind\subset \D(\g)$ denote the left
ideal in $\D(\g)$ generated by the image of
$\ind\subset\Ug$ under this homomorphism. The space
$\bigl(\D(\g)/\D(\g)\cdot\ad\ind\bigr)^{\ad\g}$ inherits
from $\D(\g)$ a natural filtered algebra structure.

One of the main results of [EG]
is a construction, for any $c\in\C$, of an  algebra
homomorphism
$$
\Phi'_c:\ \bigl(\D(\g)/\D(\g)\cdot\ad\ind\bigr)^{\ad\g}
\too \ehe.
$$
This homomorphism is compatible with the filtrations
on $\bigl(\D(\g)/\D(\g)\cdot\ad\ind\bigr)^{\ad\g}$
and $\ehe$ introduced above. So, there is a well-defined
associated graded  algebra
homomorphism 
$$\gr(\Phi'_c):\ \gr\bigl(\D(\g)/\D(\g)\cdot\ad\ind\bigr)^{\ad\g}
\too \gr(\ehe).
$$

\begin{theorem}\label{ker} For any $c\in \C$, the maps
$\Phi'_c$ and $\gr(\Phi'_c)$ are both isomorphisms.
\end{theorem}
\proof It has been explained in Remark (ii) at the end
of \cite[\S10]{EG} that both statements are immediate
consequences of Theorem \ref{t2} (the latter hasn't been
known
at the time the paper \cite{EG} was written).
\endproof

\begin{remark}
Surjectivity of the maps
$\Phi'_c$ and $\gr(\Phi'_c)$ has been already established
in \cite{EG}. 
The injectivity part was proved in  \cite{EG}
for all values of $c\in\C$ except possibly an
(unknown)  finite
set, see \cite{EG}, Theorem 7.3.
\end{remark}

\section{The functor of Hamiltonian reduction.}
\subsection{Generalities on Quantum Hamiltonian reduction.}
\label{gen}

Let $\bg$ be an arbitrary finite dimensional Lie algebra.
Given a  $\bg$-module $M$,  we write
$M^\bg:=\{m\in M\mid \rho(x)m=0,\;\forall x\in \bg\},$
for the vector space of $\bg$-invariants,
and $M_\bg:=\bg M\backslash M$, for the vector  space of
 $\bg$-{\em co}invariants.

Let $A$ be an associative
algebra, viewed as a Lie algebra
with respect to the commutator
Lie bracket. Given a Lie algebra  homomorphism
 $\rho: \bg\to A$,
one has an {\em adjoint} $\bg$-action on $A$ given by
$\ad x: a\mapsto \rho(x)\cdot a - a\cdot \rho(x),\, x\in\bg,$
$a\in A.$
The left ideal $A\cdot \rho(\bg)$ is stable under the adjoint action.
Furthermore, one shows that multiplication in $A$ induces
a well defined associative algebra structure on
$$
\A(A,\bg,\rho):=\bigl(A/A\cd\rho(\bg)\bigr)^{\ad\bg},
$$
 the space
of $\ad\bg$-invariants in $ A/ A\cd \rho(\bg)$.
The resulting algebra $\A(A,\bg,\rho)$ is called
  {\em quantum  Hamiltonian reduction} of $A$ at $\rho$.

If $A$, viewed as an $\ad\bg$-module, is semisimple, i.e.,
splits into a
(possibly infinite) direct sum of
irreducible finite dimensional $\bg$-representations, then
the operations of taking $\bg$-invariants and taking the quotient
commute, and we may write
\beq{BA}
\A(A,\bg,\rho)=\bigl(A/A\cd\rho(\bg)\bigr)^{\ad\bg}=
A^{\ad\bg}\big/(A\cd\rho(\bg))^{\ad\bg}.
\eeq
Observe that, in this formula, $(A\cdot\rho(\bg))^{\ad\bg}$ is
a {\em two-sided} ideal of the algebra $A^{\ad\bg}$.

Let $\Q:=A/ A\cdot  \rho(\bg)$, a left $A$-module.
If $a\in A$ is such that the element
$a\,\text{mod}\, A\cdot \rho(\bg) \,\in A/ A\cdot  \rho(\bg)$
is $\ad\bg$-invariant, then
the operator of right multiplication by $a$
descends to a well-defined left  $A$-linear map
$R_a:  A/ A\cdot \rho(\bg)\to A/ A\cdot \rho(\bg)$.
This gives the space $\Q$ a right
$\A(A,\bg,\rho)$-module structure, hence
makes it an $A\dash\A(A,\bg,\rho)$-bimodule.
Moreover, the right $\A(A,\bg,\rho)$-action on $\Q$
induces an algebra isomorphism
$\dis \A(A,\bg,\rho)=(\End_ A\Q)^{\opp}.
$

Let $\Lmod{A}$, resp., $\Lmod{\A(A,\bg,\rho)}$,
 be  the abelian
 category of all left $A$-modules,  resp.,
 left  $\A(A,\bg,\rho)$-modules.
Any $A$-module $M$ may be viewed also as a $\bg$-module, via the 
homomorphism $\rho$. Clearly, each of the spaces $M^\bg$ and $M_\bg$ has a natural
$A^{\ad\bg}$-module structure.

We let 
$\Lmof{(A,\bg)}$
be  the full subcategory of  $\Lmod{A}$ whose objects
are finitely generated $A$-modules which are, in addition,
completely reducible as $\bg$-modules.
We have a canonical
$A^{\ad\bg}$-module isomorphism
\begin{equation}\label{invcoinv}
M^\bg \iso M_\bg,\quad\text{for any}\en
M\in \Lmof{(A,\bg)}.\end{equation}

Further, let  $\Lmof{\A(A,\bg,\rho)}$
be the full
 subcategory of $\Lmod{\A(A,\bg,\rho)}$ whose objects are
 finitely generated $\A(A,\bg,\rho)$-modules.

\subsection{The functor $\BH$.}
We define the following functor,
called {\em Hamiltonian reduction functor}
\begin{align}\label{BH}
\BH:\ 
\Lmod{(A,\bg)}&\too\Lmod{\A(A,\bg,\rho)},\\
M&\mto \BH(M)=\Hom_{A}(\Q,M)=\Hom_{A}\bigl(A/ A\cd \rho(\bg), M\bigr)=
M^{\bg}.\nonumber
\end{align}
Here, the  action of $\A(A,\bg,\rho)$ on $\BH(M)$
comes from the tautological
{\em right} action
of $\End_ A\Q$ on $\Q,$
via the above mentioned
isomorphism $\dis \A(A,\bg,\rho)=(\End_ A\Q)^{\opp}.
$

\begin{proposition}\label{adjoint}
 Assume that $A$ is a left Noetherian algebra
and, moreover, that $A$ is a semisimple  $\ad\bg$-module. Then

\vi The algebra $\A(A,\bg,\rho)$ is left Noetherian
and $\Q$ is an object of $\Lmof{(A,\bg)}$.

\vii The   functor \eqref{BH} induces an 
{\em exact} functor
$\dis\BH: \Lmof{(A,\bg)}\to\Lmof{\A(A,\bg,\rho)}.$

\viii The functor $\ham$ in \vii has a left adjoint functor
$$
{}^\top\ham:\
\Lmof{\A(A,\bg,\rho)}\map\Lmof{(A,\bg)},\quad
E\mto \Q\otimes_{\A(A,\bg,\rho)}E.
$$ 
Furthermore, the canonical adjunction morphism 
$E\too\ham({}^\top\ham(E))$ is an isomorphism for
any $E\in \Lmof{\A(A,\bg,\rho)}.$
\end{proposition}

\proof First of all,  
we observe that the left $\bg$-action
on $\Q=A /A\cdot\rho(\bg)$ coincides with the $\ad\bg$-action
on $\Q$. 
The adjoint  action of $\bg$ on $A$, hence on $\Q$, is completely 
reducible.
We deduce that $\Q$  is   completely reducible
as a left $\bg$-module. This proves the second claim of part (i).
Also, we deduce the following natural
isomorphisms
of left $A^{\ad\bg}$-modules
 (the rightmost isomorphism below
is due to \eqref{invcoinv}):
\beq{hamq}
\A(A,\bg,\rho)=
\bigl(A /A\rho(\bg)\bigr)^{\ad\bg}=\Q^{\ad\bg}=\Q^{\bg}=\ham(\Q)\cong
\Q_\bg.
\eeq

Below, we use the notation $\A:=\A(A,\bg,\rho)$, and
 identify this  algebra   with  a quotient of the algebra $A^{\ad \bg},$
via \eqref{BA}, since the
 $\ad\bg$-action on
$A$ is  completely reducible.
Observe that the action of $A^{\ad \bg}$ on each of  the objects in
\eqref{hamq} descends to an action of  the quotient algebra $\A$. Thus, 
we may view \eqref{hamq} as a chain of isomorphisms of left
$\A$-modules.

The classical argument
due to Hilbert shows that
if $A$ is a noetherian algebra, then so is  $A^{\ad \bg}$.
A similar argument shows that,
for any $M\in \Lmof{A}$, the space
$M^\bg$ is a finitely generated  $A^{\ad \bg}$-module.
Furthermore, the functor $M\mto M^\bg$ is  exact on
the category of completely reducible $\bg$-modules.
Therefore, similar statements hold for the algebra $\A$,
a quotient of $A^{\ad \bg}.$
Parts (i)-(ii) of the Proposition follow.

Next, let $E$ be a finitely generated left $\A $-module.
Then $\Q\otimes_{\A }E$ is clearly finitely generated
over $A $. Moreover, $\Q\otimes_{\A }E$ is isomorphic,
as a left $\bg$-module, to a quotient of a direct sum
of finitely many copies of $\Q$. The latter $\bg$-module being
 completely reducible  by (i), we conclude
that $\Q\otimes_{\A }E$ is a completely reducible 
 left $\bg$-module.

Further, by general `abstract nonsense', there is  a canonical isomorphism:
$$\Hom_{A }(\Q\otimes_{\A }E, L)=
\Hom_{A }\bigl((A /A\cdot \rho(\bg))\otimes_{\A }E, L\bigr) \cong 
\Hom_{\A }(E,\, L^{\bg}),\en\forall
E\in \Lmod{\A },\,L
\in\Lmod{(A ,\bg)}.
$$
The isomorphism shows that the functor ${}^\top\ham$ is indeed a left
adjoint of $\ham$.

Using the isomorphisms of $\A$-modules from
\eqref{hamq}, we obtain
$\Tor^{\A }_1(\bg{\Q}\backslash{\Q},-)=\Tor^{\A}_1(\Q_\bg,-)=
\Tor^{\A }_1(\A ,-)=0.$
Therefore, for any
$\A $-module $E$, one has a short exact sequence
$$0\map \bg\Q\otimes_{\A }E\map
\Q\otimes_{\A }E\map(\bg{\Q}\backslash{\Q})\otimes_{\A }E
\map 0.
$$
Thus, by exactness of $\BH$ and the rightmost isomorphism in \eqref{hamq},
 we compute
\begin{align*}
\ham({}^\top\ham(E))=\left(\Q\otimes_{\A }E\right)_\bg
&=\bg(\Q\otimes_{\A }E)\backslash(\Q\otimes_{\A }E)=
(\bg\Q\otimes_{\A }E)\backslash(\Q\otimes_{\A }E)=
\\
&=
(\bg{\Q}\backslash{\Q})\otimes_{\A }E=\A \otimes_{\A }E=E.
\end{align*}
This completes the proof of part (iii) of the Proposition.
\endproof

Write $\Ker\ham$ for the full subcategory
of  $\Lmof{(A,\bg)}$ formed by the objects
$L$ such that $\ham(L)=0$. Since $\ham$ is exact, the category
$\Ker\ham$ is
a {\em Serre subcategory} in  $\Lmof{(A,\bg)}.$
Let $\Lmof{(A,\bg)}/\Ker\ham$
be the corresponding quotient category.

\begin{corollary}\label{serre} 
The functor $\ham$
  induces an equivalence
$\dis\Lmof{(A,\bg)}/\Ker\ham\iso\Lmof{\A(A,\bg,\rho)}.$
\end{corollary}
\begin{proof}
The equivalence of categories stated in the
proposition is known, by  `abstract nonsense',
 to be a formal consequence
of the existence of a left adjoint functor,
${}^\top\ham,$ such that the canonical adjunction gives
an isomorphism of functors
$\Id_{\Lmof{(A,\bg)}}\iso\ham\circ {}^\top\ham.$
The latter  isomorphism is nothing but 
Proposition \ref{adjoint}(iii).
\end{proof}

\subsection{} We return to the setting of Cherednik algebras.

 Fix $c\in\C$. We put $\bg:=\gc,\,A:=\dx,$
and let $\rho: \gc\to \dx$ be the tautological imbedding.
The algebra
$\dx$ is clearly both left and right noetherian.
The Hamiltonian reduction algebra 
$\A(\dx,\gc,\rho)=\bigl(\dx/\dx\cdot\gc\bigr)^{\ad\gc},$ c.f. \eqref{inv},
is isomorphic to the spherical Cherednik algebra
$\ehe$, by our main  Theorem \ref{iso}.

Thus, applying Proposition \ref{adjoint} and Corollary \ref{serre}
 in our present setting, we obtain the following result

\begin{proposition}\label{hamehe}
\vi The left $\dx$-module
$\Q= \dx/\dx\cdot\gc$ is an object of the
abelian category $\Lmof{(\dx,\gc)}$.

\vii The  Hamiltonian reduction functor gives an
{\em exact} functor
$$\BH:\
\Lmof{(\dx,\gc)}\too\Lmof{\ehe},\quad M\mto \Hom_{\dx}(\Q,M)=M^\g.
$$
This functor induces an equivalence
$\Lmof{(\dx,\gc)}/\Ker\BH\iso \Lmof{\ehe}.$\qed
\end{proposition}

Further, from definitions, one easily derives the following
(cf. \cite[Proposition 7.6]{EG}).

\begin{proposition}\label{Fham}
The functor $\ham$ intertwines Fourier transforms
of $\D$-modules and $\ehe$-modules, respectively, i.e., there is
a natural  isomorphism of functors
$$\ff_\hh\circ\ham\cong\ham\circ\ff_\D.\qquad\Box$$
\end{proposition}

\subsection{The Harish-Chandra $\D$-module and category ${\mathcal O}(\ehe)$.}
Recall that $\bz_+$ denotes the augmentation ideal in
$\bz\cong(\Sym\g)^{\ad\g}$,
the algebra  of $\ad\g$-invariant constant coefficients differential
operators on $\g$.
We set 
$$\SF:=\Q/\Q\cd\bz_+=\dx\big/\bigl(\dx\cd\gc+\dx\cd\bz_+).
$$

It is clear that $\SF$ is an admissible left
$\dx$-module; it may be called 
{\em Harish-Chandra}  $\dx$-module.
This name is
motivated by the works \cite{HK} and \cite{Gi2}, where the authors
considered a similar
$\D$-module for $\Ad G$-invariant eigen-distributions
on an arbitrary semisimple Lie algebra.
The  analogy with {\em loc. cit.} will be studied
further in \cite{Gi3}.

Next, let $I :=(\Sym\h)^W_+$ be the augmentation ideal
in the algebra $(\Sym\h)^W$. Given an algebra $A$ and an algebra
imbedding $ (\Sym\h)^W\into A$, we will use the notation
$\aug:=A\cdot I $ for the left ideal in $A$ generated by $I $.

Recall that the algebras
$(\Sym\h)^W$ and $\C[\h]^W$
may be viewed as subalgebras in $\ehe$.
In particular, $I $ is a subalgebra of $\ehe$.

The space 
${\ehe/\ehe\cdot(\Sym\h)^W_+}=\ehe/\aug$
has an obvious left $\ehe$-module structure.
We also consider the
following left $\hh_c$-module
$\pp:=\hh_c\otimes_{(\Sym\h)\#W}(\Sym\h/\aug)$.

\begin{lemma}\label{eM}
We have a natural $\ehe$-module isomorphism
$\ehe/\aug\cong\e\pp.$
\end{lemma}
\begin{proof}
Consider the cross product algebra $(\Sym\h)\#W$
and the natural algebra imbeddings
$(\Sym\h)^W\into (\Sym\h)\#W\into \hh_c$.
We get a  map
$f:(\Sym\h)^W/I \to\Sym\h/\Sym\h\cdot I =\Sym\h/\aug$. Tensoring this 
map with the
obvious 
inclusion $r: \ehe\into \e\hh_c$ we obtain
 a chain of  maps
\begin{align}\label{ep}
\ehe/(\ehe)\cd I =
\ehe\o_{(\Sym\h)^W}\big((\Sym\h)^W/I \big)
\stackrel{r\o f}\too 
&\e\hh_c\otimes_{(\Sym\h)\#W}(\Sym\h/\aug)\\
=&
\e\Big(\hh_c\otimes_{(\Sym\h)\#W}(\Sym\h/\aug)\Big)=
\e \pp.\nonumber
\end{align}

All the maps in this chain are filtration preserving morphisms
 of left $\ehe$-modules.
We claim that the composite morphism is, in effect, a bijection.

To prove the claim,
we consider the corresponding associated graded map.
We have $\gr(\e\hh_c)=\C[\h\times\h^*]$ is a projective, hence flat,
$\C[\h^*]\#W$-module.
We deduce, using the
identification $\C[\h^*]=\Sym\h,$
that $\e\hh_c$  is a flat right $(\Sym\h)\#W$-module, moreover,
we have
$$\gr\bigl(\e\hh_c\o_{(\Sym\h)\#W}(\Sym\h/\aug)\bigr)=
\gr(\e\hh_c)\o_{\C[\h^*]\#W}(\C[\h^*]/\aug)=
\C[\h\times\h^*]\o_{\C[\h^*]\#W}(\C[\h^*]/\aug).
$$

Similarly, we obtain $\gr(\ehe/\aug)=\C[\h\times\h^*]^W/\aug.$ Thus, 
we compute
\begin{align*}
\gr(\ehe/\aug)=\C[\h\times\h^*]^W/\aug&=(\C[\h\times\h^*]\o_W\triv)/\aug\\
&=\bigl(\C[\h\times\h^*]\o_{(\Sym\h)\#W}\Sym\h\bigr)/\langle I\rangle\\
&=
\C[\h\times\h^*]\o_{(\Sym\h)\#W}(\Sym\h/\aug)\\
&
=\gr\bigl(\e\hh_c\otimes_{(\Sym\h)\#W}(\Sym\h/\aug)\bigr)=
\gr(\e\pp).
\end{align*}
We leave to the reader to check that the composite isomorphism
above is nothing but the associated graded map corresponding
to the $\ehe$-module morphism \eqref{ep}. 
\end{proof}

Let ${\mathcal O}(\ehe)$ denote 
 category ${\mathcal O}$ for the spherical subalgebra  $\ehe$, 
see \cite{BEG}. This is a full subcategory of $\Lmof{\ehe}$
whose objects are locally nilpotent
as $I $-modules.
It is clear that $\ehe/\aug$ is an object of
${\mathcal O}(\ehe)$.

Next, recall the Hamiltonian reduction functor,
see Proposition \ref{hamehe}.
Write $\Ker\ham$ for the full subcategory
of  $\md$ formed by the objects
$L$ such that $\ham(L)=0$. 

\begin{proposition} 
The Hamiltonian reduction functor restricts to an exact functor
  $\BH: \md\to{\mathcal O}(\ehe)$. The latter functor
 induces an equivalence $\md/\Ker\ham\iso{\mathcal O}(\ehe)$.

Furthermore, we have $\ham(\SF)=\ehe/\aug$.
\end{proposition}

\proof We know that
$\Phi_c(\bz_+)=(\Sym\h)^W_+$, see \eqref{psi}.
This  immediately implies that
$\ham(L)$ is a $(\Sym\h)^W_+$-locally nilpotent 
$\ehe$-module, for any $\bz_+$-locally nilpotent 
$\dx$-module $L$. We deduce that the
 functor $\ham$ takes
 category $\md$ into ${\mathcal O}(\ehe)$.

To prove the last statement of the Proposition,
we observe that the canonical $\gc$-equivariant projection
$\Q\onto \Q^\gc$ is a morphism of right $\bz$-modules.
Hence,  we have
$(\Q\cdot\bz_+)^\gc=\Q^\gc\cdot\bz_+.$
We compute
\begin{align*}
\ham(\SF)&=\ham(\Q/\Q\cd\bz_+)=\ham(\Q)/\ham(\Q\cd\bz_+)\\
&=
\Q^\gc/(\Q\cd\bz_+)^\gc=\Q^\gc/\Q^\gc\cd\bz_+=
\ehe/\ehe\cd(\Sym\h)^W_+=\ehe/\aug,
\end{align*}
where we have used the exactness of the functor $\ham$, and
the isomorphisms in \eqref{hamq}.
\endproof

\subsection{Relation to the Hilbert scheme.}
 Let $U\sset T^*\X$ be the set formed
by the quadruples $(X,Y,i,j)\in T^*\X$
such that $V=\C[X,Y]i$, that is, such that
$i$ is a cyclic vector for the pair $(X,Y)$.
It is clear that $U$ is a $G$-stable Zariski open subset in
$T^*\X$. Furthermore, the group $G$ acts {\em freely}
on $U$ and there is a {\em universal geometric
quotient morphism}
$\Upsilon: U \stackrel{G}\too\Hilb$, 
where $\Hilb$ denotes the Hilbert scheme
of $n$ points in the plane, see \cite{Na}.

An  irreducible
component $Z$ of the Lagrangian scheme
$\La$ is said to be  {\em stable} if the set
$Z\cap U$ is dense in $Z$. 
In such a case, $Z\cap U$ is a $G$-stable closed subset
in $U$, hence, we have $Z\cap U=\Upsilon\inv(Z^{\text{Hilb}})$,
where $Z^{\text{Hilb}}:=\Upsilon(Z)$,  a closed Lagrangian
subscheme in $\Hilb.$

Now, let
$L\in \scr C_c$ and let
 $\CC(L)=\sum m_k\cdot Z_k,$ be the characteristic cycle
of $L$,
a formal integral combination of closed irreducible
subvarieties $Z_k\sset T^*\X$.
We define $\CC^{\text{Hilb}}(L):=\sum_{\{Z_k\;\text{is stable}\}}\, 
m_k\cdot Z^{\text{Hilb}}_k$, the
 formal integral combination of Lagrangian
subschemes in $\Hilb$ corresponding,
as explained above, to the
stable irreducible components $Z_k$. Thus,
$\CC^{\text{Hilb}}(L)$ is a Lagrangian cycle in $\Hilb$.

On the other hand, in the recent paper \cite{GS},
Gordon and Stafford have attached to any object
$E\in  {\mathcal O}(\ehe)$ a Lagrangian cycle 
$\CC^{GS}(E)$ in $\Hilb.$ The construction
used in \cite{GS} is totally different from the
approach of the present paper. Nevertheless,
it is likely (cf. also \cite{Gi3}) that one has:
$$\CC^{\text{Hilb}}\bigl({}^\top\ham(E)\bigr)=
\CC^{GS}(E),\quad\forall E\in {\mathcal O}(\ehe).
$$

\section{Appendix: A remark on a theorem of M. Haiman}

\centerline{\large{\sc 
Victor Ginzburg}}

\vskip 5mm

\subsection{Main result.}\label{main} Write 
 $\C[\x,\y]:=
\C[x_1,\ldots,x_n,y_1,\ldots,y_n]$
for a polynomial ring in two sets of variables
 $\x=(x_1,\ldots,x_n)$ and $\y=(y_1,\ldots,y_n)$.
The  Symmetric group $S_n$ acts naturally on the $n$-tuples
$\x$ and $\y$,
and this gives rise to an $S_n$-diagonal action
on the algebra $\C[\x,\y].$
We write $\c^{S_n}\sset \c$ for the subalgebra
of $S_n$-invariant polynomials and $A:=\c^\eps\sset\c$
for the subspace of $S_n$-alternating
polynomials. 
The space $A$
is stable under multiplication by elements of the algebra
$\c^{S_n}$, in particular, it may be viewed as a module
over $\C[\y]^{S_n}\sset\c^{S_n}$, 
 the subalgebra
of symmetric polynomials in the last $n$ variables $y_1,\ldots,y_n$.

For each $k=1,2,\ldots,$ let
$A^k$ be the $\C$-linear subspace in $\c$
spanned by the products of $k$ elements of $A$. The action of
$\C[\y]^{S_n}$
on $A$ induces one on $A^k$, hence each  space $A^k,\,k=1,2,\ldots,$
acquires a natural $\C[\y]^{S_n}$-module structure.

The goal of this Appendix is to give a direct proof of 
 the following  special case of
a much stronger result due to M. Haiman  \cite[Proposition 3.8.1]{Ha2}. 

\begin{theorem}\label{haiman} For each  $k=1,2,\ldots,$ 
the space $A^k$ is a free 
 $\C[\y]^{S_n}$-module.\end{theorem}

In  an earlier paper, Haiman showed, cf.
\cite[Proposition 2.13]{Ha1}, that  the above
theorem holds for all $k\gg0$.
The corresponding statement for {\em all}
$k$ follows from Haiman's proof 
of his
 his Polygraph theorem, the main technical result in \cite{Ha2}.

\subsection{Geometric interpretation of $A^k$.} 
The group $G=GL(V)$  acts naturally on $\M$, see \S\ref{geometry}.
Thus, we get a $G$-action $g: f\mapsto g(f),$ 
on the coordinate ring $\C[\M]$ by algebra
automorphisms.

For each $k=1,2,\ldots,$ we set
$$ \C[\M]^{(k)}:=\{f\in\C[\M]\mid g(f)=(\det g)^k \cd f,\quad\forall
g\in G\}.
$$

We will use the  notation introduced in \S\ref{state}.
There is an obvious identification
$\C[\h\times\h]=\c$.
In particular, we may view 
the vector space $A^k$, see \S\ref{main}, as a subspace in $\C[\h\times\h].$

A key ingredient in our approach to Theorem \ref{haiman} is
the following result
\begin{proposition}\label{A^k} For each $k=1,2,\ldots,$
restriction of functions  via the imbedding $\beps$, see \eqref{j},
induces a vector space isomorphism $\beps^*: \C[\M]^{(k)}\iso A^k$.
\end{proposition}

\begin{remark} Lemma \ref{f}  may be viewed as a version
of Proposition \ref{A^k} for $k=0.$
\end{remark}

\begin{remark}
For each $k=1,2,\ldots,$
M. Haiman constructed in \cite{Ha1} a natural map $A^k\to
\Gamma(\Hilb,\oo(k))$, where $\oo(1)$ is a natural
ample line bundle on $\Hilb,$ cf. \cite{Ha1}.
Moreover, it follows from the results of \cite{Ha2}
that this map is, in effect, an isomorphism.
\end{remark}

\subsection{Proof of Proposition \ref{A^k}}
 Fix nonzero volume elements
$\bv\in \wedge^n V$ and $\bv^*\in \wedge^n V^*,$ respectively.
Given an $n$-tuple $\bbf=(f_1,\ldots,f_n),\,f_r\in\C\langle x,y\rangle,$
of  noncommutative polynomials in two variables, we
consider polynomial functions $\psi,\phi\in\C[\g\times\g\times
V\times V^*]$
of the form 
\begin{align}\label{bv}
&\psi_\bbf(X,Y,i,j)=\langle\bv^*,\, f_1(X,Y)i\wedge\ldots
f_n(X,Y)i\rangle,\\
 &\phi_\bbf(X,Y,i,j)=\langle
jf_1(X,Y)
\wedge\ldots
jf_n(X,Y),\,\bv\rangle,\nonumber
\end{align}
where $f_r(X,Y)$ denotes the matrix obtained by plugging the two
matrices $X,Y\in\g$ in the  noncommutative polynomial $f(x,y)$.
We will keep the notation
$\psi_\bbf,\phi_\bbf$ for the restriction of the corresponding
function to the closed subvariety $\M\sset \g\times\g\times
V\times V^*.$ It is clear from the definition of
the imbedding $\beps:\h\times\h\into\M$ that, restricting these functions further
to the subset $\h\times\h$, one has
$\beps^*\psi_\bbf\in A$ and $\beps^*\phi_\bbf=0$. 

By Theorem \ref{t1}, we know that $\M=\M_0\cup\ldots\cup \M_n$, is a union of $n+1$
irreducible components. It is immediate from the
definition 
of the set $\M_r$, cf. \S\ref{int1}, that, for any choice of 
 $n$-tuple $\bbf=(f_1,\ldots,f_n),$ the function
$\psi_\bbf$ vanishes on $\M_r$ whenever $r\neq 0$, while
$\phi_\bbf$ vanishes on $\M_r$ whenever $r\neq n$.
Since each irredicible component is {\em reduced}, by Theorem \ref{t1},
the above vanishings hold
scheme-theoretically:
\beq{van}
\psi_\bbf|_{\M_r}=0\quad\forall r\neq0,
\quad\text{and}\quad\phi_\bbf|_{\M_r}=0\quad\forall r\neq n.
\eeq

Next, similarly to $\C[\M]^{(k)},$ for each $k\in\Z$, we introduce the space
 $\C[\g\times\g\times V\times V^*]^{(k)}$ 
of polynomial functions on $\g\times\g\times V\times V^*$
that satisfy the equation $g(f)=(\det g)^k \cdot  f,\quad\forall
g\in G.$ It is clear that 
$$\psi_\bbf\in \C[\g\times\g\times V\times V^*]^{(1)},
\quad\text{resp.},\quad
\phi_\bbf\in \C[\g\times\g\times V\times V^*]^{(-1)},
\quad\forall \bbf=(f_1,\ldots,f_n).
$$

Fix $k\in \Z$ and observe  that
$\C[\g\times\g\times V\times V^*]^{(k)}$ 
is naturally a $\C[\g\times\g\times V\times V^*]^G$-module.
Applying Weyl's fundamental theorem on $GL_n$-invariants
we deduce that this $\C[\g\times\g\times V\times V^*]^G$-module
is generated by  products  of the form
$\psi_1\cdot\ldots\cdot\psi_p\cdot\phi_1\cdot\ldots\cdot\phi_q$,
where $p-q=k$ and where each factor $\psi_r$, resp. each factor
 $\phi_r$, is of the form $\psi_\bbf$, resp., $\phi_\bbf$.

The action of $G$ on $\C[\g\times\g\times V\times V^*]$
being completely reducible, we deduce that restricting
functions from $\g\times\g\times V\times V^*$ to $\M$
yields a surjection
$\C[\g\times\g\times V\times V^*]^{(k)}\onto
\C[\M]^{(k)}$. It follows that
$\C[\M]^{(k)},$ viewed as a $\C[\M]^G$-module,
is again generated by the products
$\psi_1\cdot\ldots\cdot\psi_p\cdot\phi_1\cdot\ldots\cdot\phi_q$,
with $p-q=k$. Furthermore, from \eqref{van}, we see that 
for  $k\geq 0$
we must have $p=k\en\&\en q=0$. On the other hand, for  $k\leq 0$
we must have $p=0\en\&\en q=k$.

From now on, we assume that $k\geq 1$. Thus,
 the imbedding $\M_0\into \M$ induces a bijection
$\C[\M]^{(k)}\iso \C[\M_0]^{(k)}.$ It follows that
$\C[\M]^{(k)}$  is generated, as a $\C[\M]^G$-module, by the products
$\psi_1\cdot\ldots\cdot\psi_k$.
 Since 
 $\beps^*\psi_\bbf\in A$ for any $\bbf$, we find 
 that $\beps^*(\psi_1\cdot\ldots\cdot\psi_k)\in A^k$, hence
$\beps^*(\C[\M]^{(k)})=\beps^*(\C[\M_0]^{(k)})\sset A^k.$

To prove injectivity of the restriction map
$\beps^*: \C[\M]^{(k)}\iso A^k,$ we observe that $G\cdot
\beps(\h\times\h),$ the $G$-saturation of the image of the imbedding $\beps$,
is an irreducible variety of dimension $n^2+n=\dim \M$. Furthermore, for
any diagonal matrix $Y\in\h$ with pairwise distinct eigenvalues,
we have $\C[Y]i_o=V$. Hence, a Zariski open subset of 
$G\cdot\beps(\h\times\h)$ is contained $\M'_0$, cf. \S\ref{int1}. 
Since $\M_0=\overline{\M'_0}$ and $G\cdot\beps(\h\times\h)$ is 
irreducible, we conclude that $G\cdot\beps(\h\times\h)\sset \M_0$
and, moreover, the set $G\cdot\beps(\h\times\h)$ is Zariski dense
in $\M_0$. Thus, for any $f\in \C[\M_0]^{(k)}$ such that $\beps^*(f)=0$
we must have $f=0$. This proves injectivity of the map $\beps^*$.

We have seen above that $\C[\M]^{(k)}$  is generated, as a
 $\C[\M]^G$-module, by the
$k$-fold products
$\psi_1\cdot\ldots\cdot\psi_k$. Therefore, it suffices to prove surjectivity
 of the map $\beps^*$ for $k=1$. To prove the latter, we identify
$A=\c^\epsilon$ with $\wedge^n\C[x,y]$, the $n$-th exterior power
of the vector space $\C[x,y]$ of polynomials in 2 variables.
With this identification, the space $A$ is spanned by wedge products
of the form $f_1\wedge\ldots\wedge f_n,\,f_1,\ldots,f_n\in\C[x,y]$.

Now, by definition of the irreducible component $\M_0,$
 for any $(X,Y,i,j)\in\M_0$, we have
$[X,Y]=[X,Y]+ij=0$. Therefore, for any $f\in\C[x,y]$,
the expression $f(X,Y)$ is a well-defined matrix. In other words, for
any lift of $f$ to a noncommutative polynomial $\hat{f}\in\C\langle
x,y\rangle,$
i.e., for any $\hat{f}$ in the preimage of $f$ under the natural projection
$\C\langle
x,y\rangle\onto\C[x,y]$, we have $\hat{f}(X,Y)=f(X,Y).$
Thus, given an $n$-tuple $f_1,\ldots,f_n\in\C[x,y]$,
we have a well-defined element
$$\psi_\bbf=\langle\bv^*, f_1(X,Y)i\wedge\ldots
f_n(X,Y)i\rangle\in \C[\M]^{(1)}.
$$
It is straightforward to verify that for such an element one has
$\beps^*\psi_\bbf=f_1\wedge\ldots\wedge f_n.$
This proves 
 surjectivity of the map $\beps^*$ 
and completes the proof of the Proposition.\qed

\subsection{Proof of Theorem \ref{haiman}.} 
We have the  standard grading
$\dis\C[\y]^{S_n}=\oplus_{d\geq 0}\C^d[\y]^{S_n}$,
 by degree of the polynomial.
Write $\dis\C[\y]^{S_n}_+:=\oplus_{d> 0}\C^d[\y]^{S_n}$ 
for the augmentation ideal.

In general, let $E$ be any  flat
nonnegatively graded  $\C[\y]^{S_n}$-module.
Then, choosing representatives in $E$ of a
 $\C$-basis of the vector space
$\dis E/\C[\y]^{S_n}_+E$ yields a {\em free} $\C[\y]^{S_n}$-basis in $E$.
Hence, any  flat
nonnegatively graded  $\C[\y]^{S_n}$-module is {\em free}.

Next, we have a  $\C^\times$-action on $\M$ by dilations, given
for any $z\in\C^\times$ by the assignment
 $(X,Y,i,j)\mapsto (z\cdot X,z\cdot Y,z\cdot i,z\cdot j)$.
This  $\C^\times$-action gives rise to a natural grading 
$\dis\C[\M]=\oplus_{d\geq 0}\C^d[\M]$, on the algebra
$\C[\M]$.
With this grading, the pull-back morphism $\pi^*: \C[\y]^{S_n}
\map \C[\M]$ induced by the map $\pi: (X,Y,i,j)\mapsto \Spec Y,$
is a graded algebra morphism, so $\C[\M]$ may be viewed as
a {\em graded} $\C[\y]^{S_n}$-module.

A key point of the proof is that
the map $\pi: \M\too \C^{(n)},\,(X,Y,i,j)\mto\Spec Y$ is a
flat morphism, by  Corollary \ref{flat}.
In algebraic terms, this means
that $\C[\M]$ is a  flat (nonnegatively graded)  $\C[\y]^{S_n}$-module.
As has been explained at the beginning of the proof,
this implies that $\C[\M]$ is  free over $\C[\y]^{S_n}$.
Further, by complete reducibility of the $G$-action on $\C[\M]$,
we deduce that $\C[\M]^{(k)}$, viewed as a graded $\C[\y]^{S_n}$-submodule
in $\C[\M]$, splits off as a direct summand,
for any $k=1,2,\ldots$. Hence,
$\C[\M]^{(k)}$ is projective, in particular, flat over $\C[\y]^{S_n}$,
 as a direct summand of  a free $\C[\y]^{S_n}$-module.

To complete the proof of 
Theorem \ref{haiman} we observe
 that $\C[\M]^{(k)}$ is again a positively graded
 $\C[\y]^{S_n}$-module. Thus,  we conclude as above that 
this graded module must be free
 over $\C[\y]^{S_n}$.

\setcounter{equation}{0}
\footnotesize{

\smallskip

\noindent
{\bf W.L.G.}:
Department of Mathematics, Massachusetts Institute of Technology,
Cambridge, MA 02139, USA;\\
\hphantom{x}\quad\, {\tt wlgan@math.mit.edu}
                                                                                
\smallskip
      
\noindent                                                                          
{\bf V.G.}: 
Department of Mathematics, University of Chicago,
Chicago, IL 60637, USA;\\
\hphantom{x}\quad\, {\tt ginzburg@math.uchicago.edu}}

\end{document}